\documentstyle[11pt]{article}
\begin{document}
\newtheorem{prop}{Proposition}[section]
\newtheorem{Def}{Definition}[section]
\newtheorem{theorem}{Theorem}[section]
\newtheorem{lemma}{Lemma}[section]

\title{\bf Global well-posedness below energy space for the 1D Zakharov system}
\author{{\bf Hartmut Pecher}\\Fachbereich Mathematik\\Bergische Universit\"at -- Gesamthochschule Wuppertal\\ Gau{\ss}str. 20 \\ D-42097 Wuppertal\\ Germany\\ e-mail Hartmut.Pecher@math.uni-wuppertal.de}
\date{}
\maketitle
\begin{abstract}
The Cauchy problem for the 1-dimensional Zakharov system is shown to be globally well-posed for large data which not necessarily have finite energy. The proof combines the local well-posedness result of Ginibre, Tsutsumi, Velo and a general method introduced by Bourgain to prove a similar result for nonlinear Schr\"odinger equations.
\end{abstract}

\normalsize
\setcounter{section}{-1}
\section{Introduction}
Consider the Cauchy problem for the (1+1)-dimensional Zakharov system
\begin{eqnarray}
\label{0.1}
iu_t + u_{xx} & = & nu \\
\label{0.1a}
n_{tt} - n_{xx} & = & (|u|^2)_{xx} \\
u(0) = u_0 \, , \, n(0) & = & n_0 \, , \, n_t(0) = n_1 \label{0.2}
\end{eqnarray}
where $u$ is a complex-valued and $n$ a real-valued function defined for $(x,t) \in {\bf R}\times{\bf R}^+$.\\
The main result shows global well-posedness of the problem for rough data
$$ (u_0,n_0,n_1) \in H^{s,2}({\bf R}) \times L^2({\bf R}) \times \dot{H}^{-1,2}({\bf R}) \qquad  \mbox{with}  \qquad  1>s> 9/10 $$ without any smallness assumption.
The same result for $s=1$ is a direct consequence of the local well-posedness shown by \cite{GTV} and the conservation laws satisfied for solutions of (\ref{0.1}),(\ref{0.1a}),(\ref{0.2}), namely conservation of $ \| u(t) \| $ and
$$ E(u,n):= \|u_x(t)\|^2 + 1/2(\|n(t)\|^2 + \|V(t)\|^2) + \int^{+\infty}_{-\infty} n(t) |u(t)|^2 \, dx $$
where $ V_x = -n_t $. \\
Local well-posedness for $s > 9/10 $ follows from \cite{GTV} so that the problem is to show that the local solution exists globally in time.\\
Local and global well-posedness in dimension 2+1 and 3+1 for finite energy solutions was shown in \cite{B3}. \\
Our proof uses Bougain's ideas who introduces a general method to show global well-posedness for some types of nonlinear evolution equations for data with less regularity than needed for an application of the conservation laws directly. He applied it to the (2+1)- and (3+1)-dimensional Schr\"odinger equation
\cite{B1},\cite{B2}. Later it was also used for other model equations \cite{CST},\cite{FLP},\cite{KPV},\cite{KT}.\\
This paper is organized as follows. In section 1 the needed estimates for the nonlinearities in the $X^{s,b}$--spaces introduced by Bourgain and the $Y^s$--spaces introduced by Ginibre, Tsutsumi and Velo \cite{GTV} are given along the lines of \cite{GTV}.\\
For an equation of the form
\begin{equation}
\label{0.3}
iu_t + \phi(-i\partial_x)u = 0
\end{equation}
where $ \phi $ is a measurable real-valued function, let $ X^{s,b}$ be the completion of $ S({\bf R}^2) $ with respect to
\begin{eqnarray*}
\|f\|_{X^{s,b}} & := & \| e^{-it\phi(-i\partial_x)}f \|_{H^b_t({\bf R},H^s_x({\bf R}))} \\
& = & \| <\xi>^s <\tau>^b {\cal F}(e^{-it\phi(-i\partial_x)}f(x,t)) \|_{L^2_{x,t}} \\
& = & \| <\xi>^s <\tau+\phi(\xi)>^b \widehat{f}(\xi,\tau) \|_{L^2_{\xi,\tau}}
\end{eqnarray*}
$\dot{X}^{s,b} $ is defined similarly by replacing $<\xi>^s$ by $|\xi|^s$.\\
Similarly let $Y^s$ be the completion of $S({\bf R}^2)$ with respect to 
\begin{eqnarray*}
\|f\|_{Y^s} & := & \| <\xi>^s <\tau>^{-1} {\cal F}(e^{-it\phi(-i\partial_x)}f(x,t)) \|_{L^2_{\xi}L^1_{\tau}} \\
& = & \| <\xi>^s <\tau+\phi(\xi)>^{-1} \widehat{f}(\xi,\tau) \|_{L^2_{\xi}L^1_{\tau}}
\end{eqnarray*}
In our case we shall use these spaces for the phase functions $ \phi(\xi) = \xi^2 $ and $ \phi(\xi) = \pm |\xi| $.\\
We also have to use the norms in $ X^{s,b}(I) $ for a given time interval I defined as
$$ \| f \|_{X^{s,b}(I)} = \inf_{\tilde{f}_{|I}=f} \| \tilde{f} \|_{X^{s,b}} \quad \mbox{and similarly} \quad \| f \|_{Y^s(I)} = \inf_{\tilde{f}_{|I}=f} \| \tilde{f} \|_{Y^s} $$
In section 2 we transform the system in a standard way into a first order system for $(u,n_+,n_-)$. For details we again refer to \cite{GTV}. Using Bourgain's ideas we split the datum $ u_0 $ into a sum $ u_{01} + u_{02} $ where the low frequency part $u_{01} $ is regular and has large $ H^1 $ - norm whereas the high frequency part $u_{02} $ is just in $H^s$ with small $L^2$ - norm.\\
In section 3 the solution $(\widetilde{u},\widetilde{n_+},\widetilde{n_-})$ of the problem with data $ (u_{01},n_{0+}, n_{0-}) $ is further investigated on a suitable time interval I depending on s using the energy bounds.\\
In section 4 we consider the system fulfilled by $ (v,m_{\pm}) = (u-\widetilde{u},n_{\pm}- \widetilde{n_{\pm}}) $ with data $ (u_{02},0,0) $ and construct a solution in the same time interval I, thus we have a solution of the original problem on I. The inhomogeneous part $w(t)$ of $v(t)$ is shown to belong to $ H^{1,2}({\bf R}) $, thus is smoother than the homogeneous part $ e^{it\partial_x^2}u_{02} $ which is just in $ H^{s,2}({\bf R}) $.\\
In section 5 we show that this process can be iterated to construct a solution on any time interval $[0,T]$. What one does is to construct a solution on time intervals of equal length $|I|$. One takes as new initial data at time $|I|$ the triple $(\widetilde{u}(|I|)+w(|I|),\widetilde{n_{\pm}}(|I|)+ m_{\pm}(|I|)) $ and repeats the argument in $ [|I|,2|I|] $. Of course in each step the involved norms have to be controlled in order to be able to choose intervals of equal length.

We collect some elementary facts about the spaces $ X^{s,b} $ and $ Y^s $. If $u$ is a solution of (\ref{0.3}) with $ u(0) = f $ we have for $ b \ge 0 $ :
\begin{equation}
\label{0.4}
\| \psi_1 u \|_{X^{s,b}} \le c \| f \|_{H^s_x}
\end{equation}
If $v$ is the solution of the problem
$$ iv_t + \phi(-i\partial_x)u = F \quad , \quad v(0) = 0 $$
we have for $ b'+1 \ge b \ge 0 \ge b' > -1/2 $:
\begin{equation}
\label{0.5}
\| \psi_{\delta}v \|_{X^{s,b}} \le c \delta^{1+b'-b} \| F \|_{X^{s,b'}}
\end{equation}
and if $ b' +1 \ge b \ge 0 \ge b' $ :
\begin{equation}
\label{0.6}
\| \psi_{\delta}v \|_{X^{s,b}} \le c (\delta^{1+b'-b} \| F \|_{X^{s,b'}} + \delta^{\frac{1}{2}-b} \| F \|_{Y^s})
\end{equation}
(for a proof see \cite{GTV}, Lemma 2.1).\\
Here the cut-off function $\psi$ is in $ C^{\infty}_0({\bf R}) $ with $ supp \,\psi \subset (-2,2) $ , $ \psi \equiv 1 $ on $ [-1,1] \, ,$  $ \psi(t) = \psi(-t) $ , $ \psi(t) \ge 0 $ , $ \psi_{\delta}(t) := \psi(t/\delta) $ if $ 0<\delta\le 1 $.\\
We have $ X^{s,b}(I) \subset C^0(I,H^s({\bf R})) $ if $ b > 1/2 $ , $ I \subset {\bf R} $.\\
Moreover if $ w(t) = \int_0^t e^{i(t-s)\partial^2_x} F(s) \, ds $ and $ f(s) = e^{-is\partial^2_x} F(s) $ we have by \cite{GTV}, Lemma 2.2, especially (2.35)
\begin{eqnarray}
&&\| w(t) \|_{L^2({\bf R})}
  =  \| \int^t_0 e^{-is\partial^2_x} F(s) \, ds \|_{L^2({\bf R})}
= \| \int^t_0 f(s) \, ds \|_{L^2({\bf R})}  \\ \label{4.14a} && \le  c \| < \tau >^{-1} \widehat{f}(\tau) \|_{L^2_{\xi} L^1_{\tau}}
= c \| <\tau>^{-1} {\cal F}(e^{-is\partial_x^2} F(s))\|_{L^2_{\xi} L^1_{\tau}} =  c \| F \|_{Y^0(I)} \nonumber
\end{eqnarray}
if $ t\in [0,|I|] $ with $|I| \le 1$.\\
Thus if $ F \in Y^0(I) $ we have $ w \in C^0(I,L^2({\bf R})) $ and if $|I| \le 1$
$$ \| w \|_{L^{\infty}(I,L^2({\bf R}))} \le c \| F \|_{Y^0(I)} $$
Similarly if  $ F \in Y^1(I) $ we have $ w \in C^0(I,H^{1,2}({\bf R})) $ and if $|I| \le 1$
\begin{equation}
\label{4.14b}
\| w \|_{L^{\infty}(I,H^{1,2}({\bf R}))} \le c \| F \|_{Y^1(I)}
\end{equation}
Next we need an interpolation property for the spaces $ X^{s,b}(I) $. It is well-known that
\begin{equation}
\label{0.7}
H^b_t({\bf R},H^s_x({\bf R})) = (H^{b_0}_{t_0}({\bf R},H^{s_0}_x({\bf R})),H^{b_1}_{t_1}({\bf R},H^{s_1}_x({\bf R})))_{[\theta]}
\end{equation}
where
\begin{equation}
\label{0.8}
0 \le \theta \le 1 \quad , \quad b = (1-\theta)b_0 + \theta b_1 \quad , \quad s = (1-\theta)s_0 + \theta s_1
\end{equation}
This also holds true if ${\bf R}$ is replaced by $I$ because the restriction operator from $ H^b_t({\bf R},H^s_x({\bf R})) $ onto $ H^b_t(I,H^s_x({\bf R})) $ is a retraction with a corresponding coretraction (extension). We refer to \cite{T} here. The following interpolation property is a consequence:
\begin{equation}
\label{0.9}
X^{s,b}(I) = (X^{s_0,b_0}(I),X^{s_1,b_1}(I))_{[\theta]}
\end{equation}
with $ s,b,\theta$ as above. One only has to remark that
$$ V_{\phi} \, : \, H^b_t(I,H^s_x({\bf R})) \longrightarrow X^{s,b}(I) $$
defined by
$$ V_{\phi}f(x,t):=e^{it\phi(-i\partial_x)} f(x,t) \, , \, t \in I $$
is an isometric isomorphism.\\
Finally we have the following consequence of the Strichartz inequalities in the case of the (1+1)-dimensional Schr\"odinger equation $ \phi(\xi) = \xi ^2 $ :
\begin{equation}
\label{0.10}
\| f \|_{L^q_t({\bf R},L^r_x({\bf R}))} \le c \| f \|_{X^{0,b}}
\end{equation}
and
\begin{equation}
\| f \|_{X^{0,-b}} \le \| f \|_{L^{q'}_t({\bf R},L^{r'}_x({\bf R}))} 
\label{0.11}
\end{equation}
where $ 1/q + 1/q' = 1 \, , \, 1/r + 1/r' = 1 $ and
\begin{equation}
\label{0.12}
b_0 > 1/2 \, , \, 0 \le b \le b_0 \, , \, 1/2 \le \eta \le 1 \, , \, 2/q = 1 - \eta b/b_0 \, , \, 1/2 - 1/r = (1-\eta) b/b_0 
\end{equation}
See \cite{GTV}, Lemma 2.4 with $\nu$ = 1 , plus duality.

We use the following notation for $\lambda \in {\bf R} $:\\
$ <\lambda>:=(1+\lambda^2)^{1/2} \, , \, [\lambda]_+ := \lambda $ if $ \lambda > 0 $ , $ = \epsilon $ if $ \lambda = 0 $ , $ = 0 $ if $ \lambda < 0 $.\\
{\bf Acknowledgement:} I am grateful to A. Gr\"unrock for many helpful discussions and to T. Tao for informing the author of a gap in an earlier version of this paper.
\section{Nonlinear estimates}
Our aim here is to estimate the nonlinearities $ f = n_{\pm} u $ in $ X^{s,-a_1} $ for given $ n_{\pm} \in X^{l,a} $ and $ u \in X^{k,a_2} $ for suitable $ s,l,k,a_1,a,a_2$ and also in $ Y^s $. We estimate $ \widehat{f}(\xi'_1,\tau_1) = (\widehat{n_{\pm}} \ast \widehat{u})(\xi'_1,\tau_1) $ in terms of $ \widehat{n_{\pm}}(\xi,\tau) $ and $ \hat{u}(\xi'_2,\tau_2), $ where $ \xi = \xi'_1 - \xi'_2 ,$ $ \tau = \tau_1 - \tau_2 $. We also introduce the variables  $ \sigma_1 = \tau_1 + {\xi'_1}^2 $ , $ \sigma_2 = \tau_2 + {\xi'_2}^2 $ , $ \sigma = \tau \pm |\xi| $ so that
\begin{equation}
\label{1.1}
z := {\xi'_1}^2 - {\xi'_2}^2\mp|\xi| = \sigma_1 - \sigma_2 -\sigma
\end{equation}
Define $ \widehat{v_2} = < \xi'_2 >^k <\sigma_2>^{a_2} \hat{u} , \hat{v} = < \xi >^l <\sigma>^a \widehat{n_{\pm}} $ so that $ \|u\|_{X^{k,a_2}} = \|v_2\|_2 $ , $ \|n_{\pm}\|_{X^{l,a}} = \|v\|_2 $ .\\
In order to estimate $f$ in $X^{s,-a_1}$ we take its scalar product with a function in $X^{-s,a_1}$ with Fourier transform $<\xi'_1>^s<\sigma_1>^{-a_1} \widehat{v_1} $ , $ v_1 \in L^2 $. In the sequel we want to show an estimate of the type
\begin{equation}
|S| \le c \|v\|_2\|v_1\|_2\|v_2\|_2
\label{1.2}
\end{equation}
where
\begin{equation}
\label{1.3}
S :=\int \frac{|\widehat{v}\widehat{v_1}\widehat{v_2}| <\xi'_1>^s}{<\sigma>^a<\sigma_1>^{a_1}<\sigma_2>^{a_2}<\xi'_2>^k<\xi>^l} \, d\xi'_1d\xi'_2d\tau_1d\tau_2
\end{equation}
This directly gives the desired estimate
\begin{equation}
\label{1.4}
\| n_{\pm}u\|_{X^{s,-a_1}} \le c \| n_{\pm} \|_{X^{l,a}} \|u\|_{X^{k,a_2}}
\end{equation}
\begin{prop}
\label{Proposition 1.1}
The estimate (\ref{1.4}) holds under the following conditions:
\begin{eqnarray}
&& k,l \ge 0 
\label{1.5}\\
\label{1.6}
&& s-k < \min(2a-\frac{1}{2},2a_1-\frac{1}{2},2(a+a_1)-\frac{3}{2}) \\
\label{1.7}
&& s-l \le 2 a_1 \\
\label{1.8}
&& a,a_1,a_2 > 1/4 \, , \, a_1 \le 1/2 \\
\label{1.9}
&& a+a_1 \, , \, a+a_2 \, , \,  a_1 + a_2 > 3/4 \\
\label{1.10}
&& k+a_1 \, , \, k+a_2 > 1/2 \quad , \quad k+a_1+a_2 > 1
\end{eqnarray}         
\end{prop}
{\bf Remark 1:} We simplify (\ref{1.1}) in the following way:\\
if (\ref{1.1}) holds with the minus-sign and if $ \xi'_1 \ge \xi'_2 $ (resp. $ \xi'_1 \le \xi'_2 $) we have
$$ z = {\xi'_2}^2-{\xi'_2}^2-|\xi'_1-\xi'_2| = (\xi'_1 \mp \frac{1}{2})^2 - (\xi'_2 \mp \frac{1}{2})^2 = \xi_1^2 - \xi_2^2 $$
where $ \xi_i := \xi'_i \mp 1/2 $.\\
Thus the region $ \xi'_1 \ge \xi'_2 $ (resp. $ \xi'_1 \le \xi'_2 $) of $S$ is majorized by
\begin{equation}
\label{1.11}
\bar{S} = \int \frac{|\widehat{v}(\xi,\tau) \widehat{v_1}(\xi_1\pm\frac{1}{2},\tau_1) \widehat{v_2}(\xi_2\pm\frac{1}{2},\tau_2)| <\xi_1>^s}{<\sigma>^a<\sigma_1>^{a_1}<\sigma_2>^{a_2}<\xi_2>^k<\xi>^l} \, d\xi_1d\xi_2d\tau_1d\tau_2
\end{equation}
where now
\begin{eqnarray}
\label{1.12}
&&z = \xi_1^2 - \xi_2^2 = \sigma_1 - \sigma_2 -\sigma \, , \, \xi = \xi_1 - \xi_2 \, , \, \tau = \tau_1 - \tau_2 \\
&&\sigma_i = \tau_i + (\xi_i\pm1/2)^2 \, ,  \, \sigma = \tau + |\xi| = \tau + |\xi_1 - \xi_2|
\nonumber
\end{eqnarray}
Also the plus-sign in (\ref{1.1}) can be treated similarly by again defining $ \xi_i = \xi'_i \pm 1/2 $. If one wants to estimate $ \bar{S} $ by $ \|v\|_2\|v_1\|_2\|v_2\|_2 $ the variables $\xi_i$ and $\xi_i \pm 1/2 $ of $\widehat{v_i}$ are completely equivalent, thus we do not distinguish between them.\\
{\bf Remark 2:} We use the following application of Schwarz' inequality: in order to estimate
$$ I = \int |\widehat{v}(\zeta)\widehat{v_1}(\zeta_1)\widehat{v_2}(\zeta_2)K(\zeta_1,\zeta_2)| \, d\zeta_1d\zeta_2 $$
where $ \zeta = (\xi,\tau) \, , \, \zeta_i = (\xi_i,\tau_i) \, , \, \zeta = \zeta_1 - \zeta_2 $ one has
\begin{eqnarray}
\nonumber
|I|^2 & \le & \|v\|_2^2 \int\left(\int|\widehat{v_1}(\zeta+\zeta_2)\widehat{v_2}(\zeta_2)K(\zeta+\zeta_2,\zeta_2)|\, d\zeta_2\right)^2 \, d\zeta \\
\label{1.13}
& \le & \|v\|_2^2 \left( \sup_{\zeta} \int |K(\zeta+\zeta_2,\zeta_2)|^2 \, d\zeta_2 \right) \int |\widehat{v_1}(\zeta+\zeta_2)\widehat{v_2}(\zeta_2)|^2 \, d\zeta_2 d\zeta \quad \\ 
\nonumber
& = & C^2 \|v\|_2^2 \|v_1\|_2^2 \|v_2\|_2^2
\end{eqnarray}
with
\begin{equation}
\label{1.14}
C^2 = \sup_{\zeta} \int |K(\zeta_1,\zeta_2)|^2 \, d\zeta_2
\end{equation}
where the integral runs over $\zeta_2$ (or $\zeta_1$) for fixed $\zeta$.\\
Similar estimates hold by circularly permuting the variables $\zeta,\zeta_1,\zeta_2$.\\
{\bf Proof of Prop. \ref{Proposition 1.1}:} We estimate (\ref{1.11}) in several subregions.\\
{\bf Region a:} $ |\xi_1| \ge 2 |\xi_2| $ \\
{\bf Case aa:} $ \sigma_1 $ dominant, i.e. $|\sigma_1|\ge |\sigma| $ , $ |\sigma_1| \ge |\sigma_2| $.\\
We show according to (\ref{1.13}),(\ref{1.14}): 
\begin{eqnarray*}
 C_1^2 :=  \sup_{\xi_1,\sigma_1} <\sigma_1>^{-2a_1}<\xi_1>^{2s} \hspace{-0.5cm}\int\limits_{\xi_1,\sigma_1 \mbox{\scriptsize{fixed}}} \hspace{-0.7cm} <\sigma>^{-2a}<\sigma_2>^{-2a_2}<\xi_2>^{-2k}<\xi>^{-2l} d\xi_2 d\sigma_2  \\
 < \infty \hspace{11.7cm}
\end{eqnarray*}
We have
\begin{eqnarray*}
\int <\sigma>^{-2a}<\sigma_2>^{-2a_2}\, d\sigma_2 & = & \int <\sigma_1-\sigma_2-\xi_1^2+\xi_2^2>^{-2a} <\sigma_2>^{-2a_2} \, d\sigma_2 \\
& \le & c <\xi_1^2-\xi_2^2-\sigma_1>^{-\alpha_1}
\end{eqnarray*}
by (\ref{1.12}) and \cite{GTV}, Lemma 4.2 with $\alpha_1:=2\min(a,a_2)-[1-2\max(a,a_2)]_+ $ if $ a+a_2 > 1/2 $ which holds by (\ref{1.8}).\\
Thus 
$$ C_1^2 \le c \sup_{\xi_1,\sigma_1}<\sigma_1>^{-2a_1}<\xi_1>^{2s} \int <\xi_2>^{-2k}<\xi_1^2-\xi_2^2-\sigma_1>^{-\alpha_1}<\xi>^{-2l}\, d\xi_2 $$
Now using $ | \xi| \ge | \xi_1|-| \xi_2| \ge \frac{1}{2}|\xi_1| $ , thus $ <\xi>^{-2l} \le c <\xi_1>^{-2l} $ , and substituting $ y = \xi_2^2 $ , $ dy = 2|y|^{1/2}d\xi_2 $ , $ |y| \le \frac{1}{4}\xi_1^2 $ we get 
\begin{eqnarray*}  C_1^2 &\le & c \sup_{\xi_1,\sigma_1} <\sigma_1>^{-2a_1}<\xi_1>^{2(s-l)} \int_{-\infty}^{+\infty} |y|^{-1/2} <y>^{-k} \\ 
&& \hspace{5cm} \chi_{\{|y|\le \frac{1}{4}\xi_1^2\}}<y-(\xi_1^2-\sigma_1)>^{-\alpha_1} \, dy 
\end{eqnarray*}
According to \cite{GTV}, Lemma 4.1 the integral takes its maximum at $ \xi_1^2 = \sigma_1 $ , so that
$$ C_1^2 \le c \sup_{\xi_1,\sigma_1} <\sigma_1>^{-2a_1} <\xi_1>^{2(s-l)} \int_{-\infty}^{+\infty}|y|^{-1/2} <y>^{-(k+\alpha_1)} \, dy $$
The integral converges if $ k+\alpha_1 > 1/2 $. By definition $ \alpha_1 = 2a $ or $ 2a_2 $ or $ 2(a+a_2)-1 $ up to an $\epsilon$-term. Now $ k+2a > 1/2 $ by (\ref{1.5}),(\ref{1.8}), similarly $ k+2a_2 > 1/2 ,$  and $ k+2(a+a_2)-1 > 1/2 $ by (\ref{1.5}),(\ref{1.9}), thus $ k+\alpha_1 > 1/2 $. Moreover $ \frac{3}{4}\xi_1^2 \le \xi_1^2-\xi_2^2 = \sigma_1 - \sigma_2 -\sigma \le 3|\sigma_1| $ so that
$$ C_1^2 \le c \sup_{\xi_1} <\xi_1>^{-4a_1+2(s-l)} < \infty $$
because $ -4a_1+2(s-l) \le 0 $ by (\ref{1.7}).\\
{\bf Case ab:} $\sigma_2$ dominant, i.e. $|\sigma_2| \ge |\sigma| $ , $ |\sigma_2| \ge |\sigma_1| $.\\
By (\ref{1.13}),(\ref{1.14}) we have to show 
\begin{eqnarray*}
C_2^2 := \sup_{\xi_2,\sigma_2} <\sigma_2>^{-2a_2} <\xi_2>^{-2k} \hspace{-0.7cm}\int\limits_{\sigma_2,\xi_2 \,\mbox{\scriptsize{fixed}}}\hspace{-0,7cm}<\xi_1>^{2s}<\sigma>^{-2a}<\sigma_1>^{-2a_1}<\xi>^{-2l}  d\xi_1 d\sigma_1 
\\
< \infty \hspace{11.5cm}
\end{eqnarray*}
We have
$$ C_2^2 \le c \sup_{\xi_2,\sigma_2} <\sigma_2>^{-2a_2} <\xi_2>^{-2k} \hspace{-0.6cm}\int\limits_{\sigma_2,\xi_2 \,\mbox{\scriptsize{fixed}}} \hspace{-0,6cm}<\xi_1>^{2(s-l)} <\sigma>^{-2a} <\sigma_1>^{-2a_1} \, d\xi_1 d\sigma_1 $$
Substituting $ \xi_1 $ by $z$ with fixed $\xi_2$ gives $ z=\xi_1^2-\xi_2^2 $ , $ dz = 2\xi_1d\xi_1 $ and by (\ref{1.12}) $ 3\xi_2^2 \le z = \sigma_1 - \sigma_2 -\sigma \le 3|\sigma_2| $. Thus
\begin{eqnarray*}
C_2^2 & \le c & \sup_{\xi_2,\sigma_2} <\sigma_2>^{-2a_2} <\xi_2>^{-2k} \\ &&\int_{3\xi_2^2}^{3|\sigma_2|} <z>^{s-l} |z|^{-1/2} \left(\int<\sigma_1-(\sigma_2+z)>^{-2a}<\sigma_1>^{-2a_1}d\sigma_1\right)dz
\end{eqnarray*}
Now the inner integral is estimated using \cite{GTV}, Lemma 4.2 by $c<\sigma_2+z>^{-\alpha_2}$, if $a+a_1>1/2$ (which follows from (\ref{1.8})) with $\alpha_2:=2\min(a,a_1)-[1-2\max(a,a_1)]_+$. This gives
\begin{eqnarray*}
C_2^2 & \le & c \sup_{\xi_2,\sigma_2} <\sigma_2>^{-2a_2}<\xi_2>^{-2k} \int_{3\xi_2^2}^{3|\sigma_2|}<z>^{s-l}|z|^{-1/2}<z+\sigma_2>^{-\alpha_2}dz \\
& \le & c \sup_{\sigma_2} <\sigma_2>^{-2a_2}\int_0^{3|\sigma_2|}<z>^{s-l}|z|^{-1/2}<z+\sigma_2>^{-\alpha_2}dz
\end{eqnarray*}
We split up the integral into the parts $0\le z \le\frac{1}{2}|\sigma_2| $ and $\frac{1}{2}|\sigma_2| \le z \le 3|\sigma_2|$. Assuming w.l.o.g. $s-l\ge 0$ we estimate the first part by $c<\sigma_2>^{s-l+\frac{1}{2}-\alpha_2}$ and the second part by $c<\sigma_2>^{s-l-\frac{1}{2}+[1-\alpha_2]_+}$ which is the larger one. Thus
$$ C_2^2 \le c \sup_{\sigma_2} <\sigma_2>^{-2a_2+s-l-\frac{1}{2}+[1-\alpha_2]_+} < \infty $$
if
\begin{equation}
\label{1.15}
s-l \le 2a_2+\frac{1}{2} - [1-\alpha_2]_+
\end{equation}
Now $ \alpha_2 = 2a $ or $ 2a_1 $ or $ 2(a+a_1)-1 $ and we have
\begin{eqnarray*}
s-l &\hspace{-0.2cm} \le &\hspace{-0.2cm} 2a_1 < 2a_1+2a_2-\frac{1}{2} \quad\mbox{by (\ref{1.7}),(\ref{1.8})} \\
s-l & \hspace{-0.2cm} \le &\hspace{-0.2cm}  2a_1 < 2a_1+[2(a+a_2)-\frac{3}{2}] = 2a_2+\frac{1}{2}+2(a+a_1)-2 \;\;\mbox{by (\ref{1.7}),(\ref{1.8}),(\ref{1.9})} \\
s-l & \hspace{-0.2cm} \le & \hspace{-0.2cm} 1 < 1+[2(a_2+a)-\frac{3}{2}] = 2a_2+2a-\frac{1}{2} \quad \mbox{by (\ref{1.7}),(\ref{1.8}),(\ref{1.9})} \\
s-l & \hspace{-0.2cm} \le &\hspace{-0.2cm}  1 < 2a_2+\frac{1}{2}\quad \mbox{by (\ref{1.7}),(\ref{1.8})}
\end{eqnarray*}
Thus (\ref{1.15}) is satisfied. \\
{\bf Case ac:} $\sigma$ dominant, i.e. $|\sigma|\ge|\sigma_1| $ , $ |\sigma| \ge |\sigma_2| $ \\
We have to show that 
\begin{eqnarray*} 
C^2 := \sup_{\xi,\sigma} <\sigma>^{-2a}<\xi>^{-2l}
 \hspace{-0.5cm}
\int\limits_{\sigma,\xi\,\mbox{\scriptsize{fixed}}}
 \hspace{-0.5cm}
<\xi_2>^{-2k}<\sigma_1>^{-2a_1}<\sigma_2>^{-2a_2}<\xi_1>^{2s} d\xi_2 d\sigma_2 \\ 
< \infty 
\end{eqnarray*}
Using $ |\xi| \ge \frac{1}{2}|\xi_1| $ and $ |\xi| \le |\xi_1|+|\xi_2| \le \frac{3}{2}|\xi_1| $, thus 
$$ \frac{1}{3} |\xi|^2 \le \frac{3}{4} \xi_1^2 \le \xi_1^2 - \xi_2^2 = z = \sigma_1 - \sigma_2 - \sigma \le 3|\sigma| $$
we get
$$ C^2 \le c \sup_{\xi,\sigma}<\xi>^{-4a+2(s-l)} \int\limits_{\sigma,\xi\,\mbox{\scriptsize{fixed}}} <\xi_2>^{-2k}<\sigma_1>^{-2a_1}<\sigma_2>^{-2a_2} d\xi_2 d\sigma_2 $$
Substituting $\xi_2$ by $z$ for fixed $\xi$ gives 
$z = \xi_1^2-\xi_2^2 = (\xi_1+\xi_2)\xi = (\xi+2\xi_2)\xi$ , $ \frac{dz}{d\xi_2} = 2\xi $ and 
$ z-\xi^2 = (\xi_1^2-\xi_2^2)-\xi^2 = (\xi_1+\xi_2)\xi-(\xi_1-\xi_2)\xi = 2\xi_2\xi $ thus $ \xi_2 = \frac{z-\xi^2}{2\xi} $ leads to
$$ C^2 \le c \sup_{\xi,\sigma} <\xi>^{-4a+2(s-l)} |\xi|^{-1} \hspace{-0.2cm}\int\limits_0^{3\xi^2} \hspace{-0.2cm}< \frac{z-\xi^2}{2\xi}>^{-2k}(\int\hspace{-0.2cm}<\sigma_1>^{-2a_1}<\sigma_2>^{-2a_2}d\sigma_2)dz $$
We used here
\begin{equation}
\label{1.17}
z = \xi_1^2 - \xi_2^2 =\xi(\xi_1 + \xi_2) \le \frac{3}{2} |\xi|\,|\xi_1| \le 3 \xi^2
\end{equation}
Now
\begin{eqnarray*}
\int <\sigma_1>^{-2a_1} <\sigma_2>^{-2a_2} d\sigma_2 & = & \int <\sigma_2-(-z-\sigma)>^{-2a_1}<\sigma_2>^{-2a_2} d\sigma_2 \\
& \le & c<z+\sigma>^{-\alpha}
\end{eqnarray*}
by \cite{GTV}, Lemma 4.2 with $\alpha=2\min(a_1,a_2)-[1-2\max(a_1,a_2)]_+ $ using $ a_1+a_2 > 1/2 $ which holds by (\ref{1.8}).\\
Substitute $ y = z - \xi^2 $ and use $|y| \le |z|+\xi^2\le4\xi^2$ by (\ref{1.17}) to conclude
\begin{eqnarray*}
C^2 & \le & c \sup_{\xi,\sigma} <\xi>^{-4a+2(s-l)} |\xi|^{-1} \hspace{-0.1cm}\int_{-\infty}^{\infty} \hspace{-0.2cm}\chi_{\{|y|\le4\xi^2\}} <\frac{y}{2|\xi|}>^{-2k} <y+\xi^2+\sigma>^{-\alpha} dy \\
& = & c \sup_{\xi} <\xi>^{-4a+2(s-l)} |\xi|^{-1} \int_0^{4\xi^2}  <\frac{y}{2|\xi|}>^{-2k} <y>^{-\alpha} dy
\end{eqnarray*}
by \cite{GTV}, Lemma 4.1.\\
If $|\xi|\le 1$ we directly get
$$ C^2 \le c \sup_{|\xi|\le 1} |\xi|^{-1} \int_0^{4\xi^2} <\frac{y}{2|\xi|}>^{-2k} dy \le c \sup_{|\xi|\le 1} |\xi|^{-1}|\xi|^2 < \infty $$
If $|\xi|\ge 1 $ we get
\begin{eqnarray*}
C^2 & \le & c \sup_{|\xi|\ge 1} |\xi|^{-4a+2(s-l)-1} \left(\int_0^{|\xi|} <y>^{-\alpha} dy + \int_{|\xi|}^{4\xi^2} |y|^{-2k-\alpha}dy\, |\xi|^{2k} \right) \\
& \le & c \sup_{|\xi|\ge 1}(|\xi|^{-4a+2(s-l)-1+[1-\alpha]_+} + |\xi|^{-4a+2(s-l)-1+2k} |\xi|^{-2k-\alpha+1}) \\
& \le & c \sup_{|\xi|\ge 1}|\xi|^{-4a+2(s-l)-1+[1-\alpha]_+} < \infty
\end{eqnarray*}
if
\begin{equation}
\label{1.18}
2k+\alpha > 1
\end{equation}
and
\begin{equation}
\label{1.19}
-4a+2(s-l)-1+[1-\alpha]_+ \le 0
\end{equation}
By the definition of $\alpha$ (\ref{1.18}) holds if $2k+2a_1>1$ and $2k+2a_2>1$ and $ 2k+2(a_1+a_2)-1>1$ which follows from (\ref{1.10}).\\
In order to show (\ref{1.19}) we use $\alpha=2a_1$ or $2a_2$ or $2(a_1+a_2)-1$ and get $ 2(s-l) \le 2 < 4a+1 $ and $ 2(s-l) \le 2 < 2(a+a_1)+2a=4a+2a_1$ and $2(s-l)\le 2 < 2(a+a_2)+2a = 4a+2a_2$ and $ 2(s-l) \le 2 < 2(a+a_1)+2(a+a_2)-1$ by (\ref{1.7}),(\ref{1.8}),(\ref{1.9}), which implies (\ref{1.19}). \\
{\bf Region b:} $|\xi_1| \le 2|\xi_2| $ \\
We show
\begin{eqnarray*}
C_2^2 := \sup_{\xi_2,\sigma_2} <\sigma_2>^{-2a_2} <\xi_2>^{-2k} \hspace{-0.7cm}\int\limits_{\sigma_2,\xi_2 \,\mbox{\scriptsize{fixed}}}\hspace{-0,7cm}<\xi_1>^{2s}<\sigma>^{-2a}<\sigma_1>^{-2a_1}<\xi>^{-2l}  d\xi_1 d\sigma_1 
\\
< \infty \hspace{11.5cm}
\end{eqnarray*}
We have
$$ C_2^2 \le c \sup_{\xi_2,\sigma_2} <\xi_2>^{-2k} \int \limits_{\sigma_2,\xi_2 \,\mbox{\scriptsize{fixed}}} <\xi_1>^{2s}<\sigma>^{-2a}<\sigma_1>^{-2a_1} d\xi_1 d\sigma_1 $$
Now by (\ref{1.12}) and \cite{GTV}, Lemma 4.2:
\begin{eqnarray*}&&\int <\sigma>^{-2a} <\sigma_1>^{-2a_1} d\sigma_1 \\&&  =  \int <\sigma_1-(\sigma_2+\xi_1^2-\xi_2^2)>^{-2a} <\sigma_1>^{-2a_1} d\sigma_1 \\&&
\le  c <\sigma_2+\xi_1^2-\xi_2^2>^{-{\alpha}_2}
\end{eqnarray*}
with $ \alpha_2 $ as above.\\ The substitution $ y = \xi_1^2 , $ $ dy = 2|y|^{1/2}d\xi_1 , $ $ y \le 4\xi_2^2 $ gives \begin{eqnarray*}
C_2^2 & \le & c \sup_{\xi_2,\sigma_2} \int <\xi_1>^{2(s-k)} <\sigma_2+\xi_1^2-\xi_2^2>^{-{\alpha}_2} d\xi_1 \\
& \le & c \sup_{\xi_2,\sigma_2} \int_0^{4\xi_2^2} <y>^{s-k} <y-(\xi_2^2-\sigma_2)>^{-{\alpha}_2} |y|^{-1/2} dy \\
& \le & c \sup_{\xi_2,\sigma_2} \int_{-\infty}^{+\infty}<y>^{s-k}|y|^{-1/2} \chi_{\{|y|\le 4\xi_2^2\}}<y-(\xi_2^2-\sigma_2)>^{-{\alpha}_2} dy \\
& \le & c \int_0^{\infty}<y>^{s-k}|y|^{-1/2} <y>^{-{\alpha}_2} dy
\end{eqnarray*}
by use of \cite{GTV}, Lemma 4.1 (remark that $s-k \le 1/2$ by (\ref{1.6}),(\ref{1.8})).\\
Thus  $ C_2^2 < \infty $  provided  $s-k <  \alpha_2- \frac{1}{2}$
which follows from (\ref{1.6}) and completes the proof.\\

Our next aim is to give a similar estimate for $ f = n_{\pm} u $ in $ Y^s $. We first integrate $<\sigma_1>^{-1}\widehat{f}$ over $\tau_1$ and take the scalar product with a function in $H_x^{-k}({\bf R})$ with Fourier transform $<\xi_1>^k\widehat{w_1}$ , $w_1 \in L_x^2({\bf R}) $. We show that an estimate of the type
$$ |\tilde{S}| \le c \|v\|_2\|w_1\|_2\|v_2\|_2 $$
holds, where
$$ \tilde{S} := \int\frac{|\widehat{v}\widehat{w_1}\widehat{v_2}|<\xi'_1>^s}{<\sigma>^a<\sigma_1><\sigma_2>^{a_2}<\xi'_2>^k<\xi>^l} d\xi'_1d\xi'_2\tau_1d\tau_2 $$
and the notation is the same as for $S$ before.\\
This directly gives the estimates
\begin{equation}
\label{1.21}
\|n_{\pm}u\|_{Y^s} \le c\|n_{\pm}\|_{X^{l,a}} \|u\|_{X^{k,a_2}}
\end{equation}
\begin{prop}
\label{Proposition1.2}
The estimate (\ref{1.21}) holds under the following conditions:
\begin{eqnarray}
\label{1.22}
k>0 \, , \, l \ge 0 \\
\label{1.23}
s-k < \min(2a-\frac{1}{2},\frac{1}{2}) \\
\label{1.24}
s-l \le 1 \\
\label{1.25}
a,a_2 > 1/4 \\
\label{1.26}
a+a_2 > 3/4 \\
\label{1.27}
k+a_2 > 1/2
\end{eqnarray}
\end{prop}
{\bf Remark:} The same remarks as for Prop. \ref{Proposition 1.1} apply. Thus we estimate
\begin{equation}
\bar{\tilde{S}} = \int \frac{|\widehat{v}(\xi,\tau)\widehat{w_1}(\xi_1\pm 1/2)\widehat{v_2}(\xi_2\pm 1/2,\tau_2)|<\xi_1>^s}{
<\sigma>^a<\sigma_1><\sigma_2>^{a_2}<\xi_2>^k<\xi>^l} d\xi_1d\xi_2d\tau_1d\tau_2
\label{1.28}
\end{equation}
where again (\ref{1.12}) holds with the same notation again.\\
{\bf Proof of Prop.1.2:} The proof works along the same lines as the foregoing one. Choose $a_1$ such that
\begin{equation}
\label{1.29}
1/2 > a_1 > \min(\frac{1}{4},\frac{3}{4}-a,\frac{3}{4}-a_2,\frac{5}{4}-a-a_2,\frac{1}{2}-k,1-a_2-k,\frac{s-k}{2}+\frac{1}{4},\frac{s-k}{2}-a+\frac{3}{4})
\end{equation}
This and (\ref{1.22})-(\ref{1.27}) imply the conditions (\ref{1.5}),(\ref{1.6}),(\ref{1.8}),(\ref{1.9}),(\ref{1.10}) and
\begin{equation}
\label{1.30}
a+a_1+a_2 > 5/4
\end{equation}
We consider the regions of integration similarly as in the proof of Prop.\ref{Proposition 1.1} with the minor variation that {\bf case aa} means $|\sigma_1|\ge4|\sigma_2|$ , $|\sigma_1|\ge4|\sigma|$ , {\bf case ab} is given by $|\sigma_2|\ge\frac{1}{4}|\sigma_1|$ , $|\sigma_2|\ge|\sigma|$ , and {\bf case ac} by $|\sigma|\ge\frac{1}{4}|\sigma_1|$ , $|\sigma|\ge|\sigma_2|$.\\
{\bf Case aa:} $|\xi_1|\ge2|\xi_2|$ \\
Choose
$$ \widehat{v_1}(\xi_1\pm1/2,\tau_1) := <\sigma_1>^{-1/2}\widehat{w_1}(\xi_1\pm1/2)\chi_{\{\frac{1}{2}|\sigma_1|\le\xi_1^2\le2|\sigma_1|\}} $$
Now in the considered region we have $ z=\xi_1^2-\xi_2^2\ge\frac{3}{4}\xi_1^2 $ and
$$ z=\sigma_1-\sigma_2-\sigma=|\sigma_1-\sigma_2-\sigma|\ge|\sigma_1|-|\sigma_2|-|\sigma|\ge|\sigma_1|-\frac{1}{4}|\sigma_1|-\frac{1}{4}|\sigma_1|=\frac{1}{2}|\sigma_1| $$
Moreover
$$ \frac{4}{3}z=\frac{4}{3}(\sigma_1-\sigma_2-\sigma)\le\frac{4}{3}(|\sigma_1|+|\sigma_2|+|\sigma|)\le\frac{4}{3}(|\sigma_1|+\frac{1}{4}|\sigma_1|+\frac{1}{4}|\sigma_1|)=2|\sigma_1|$$
altogether
$$ \frac{1}{2}|\sigma_1| \le z \le \xi_1^2 \le \frac{4}{3}z\le 2|\sigma_1| $$
This means that in the case at hand $ \chi \equiv 1 $ , so that
$$\bar{\tilde{S}} = \int \frac{|\widehat{v}(\xi,\tau)\widehat{v_1}(\xi_1\pm 1/2,\tau_1)\widehat{v_2}(\xi_2\pm 1/2,\tau_2)|<\xi_1>^s}{
<\sigma>^a<\sigma_1>^{1/2}<\sigma_2>^{a_2}<\xi_2>^k<\xi>^l} d\xi_1d\xi_2d\tau_1d\tau_2$$
where the integral runs over the region aa. This is exactly the term treated in case aa of Prop. \ref{Proposition 1.1} with $a_1$ replaced by $1/2$. Since all the assumptions of  Prop. \ref{Proposition 1.1} are satisfied with $a_1=1/2$ under our assumptions (\ref{1.22}) - (\ref{1.27}) we get from that proof an estimate by $ c\|v\|_2\|v_1\|_2\|v_2\|_2 $. Because we want to have an estimate by $ c\|v\|_2\|w_1\|_2\|v_2\|_2 $, the only thing to be checked is $\|v_1\|_2\le c\|w_1\|_2$. But this is true, namely
\begin{eqnarray*}
\|v_1\|_2^2 & = & \int\hspace{-0.9cm}\int\limits_{\frac{1}{2}|\sigma_1|\le\xi_1^2\le2|\sigma_1|} <\tau_1+(\xi_1\pm1/2)^2>^{-1}|\widehat{w_1}(\xi_1\pm1/2)|^2d\tau_1d\xi_1 \\
& = & \int(\int\limits_{\frac{1}{2}\xi_1^2\le|\sigma_1|\le2\xi_1^2} <\sigma_1>^{-1}d\sigma_1)|\widehat{w_1}(\xi_1\pm1/2)|^2d\xi_1
\end{eqnarray*}
Now the inner integral is bounded independently of $\xi_1$ by some logarithm. Thus $ \|v_1\|_2 \le c\|w_1\|_2 $. This concludes the proof of case aa.

In all other cases we define $ \widehat{v_1}:=<\sigma_1>^{a_1-1}\widehat{w_1} $ with $ a_1 $ as above. Then we simply have
$$ \|v_1\|_2^2 = \int\int<\sigma_1>^{2(a_1-1)}|\widehat{w_1}(\xi_1)|^2d\sigma_1d\xi_1 \le c\|w_1\|_2^2 $$
because $ a_1<1/2$ by our choice (\ref{1.29}).\\
Now with this choice $\bar{\tilde{S}} $ reduces to $ \bar{S} $ so that it only remains to check that the old estimates of Prop. \ref{Proposition 1.1} in all other cases (with the {\bf modified} cases aa,ab,ac as above) remain true. According to the remarks at the beginning of this proof we have to avoid using (\ref{1.7}) in its strong form, namely (\ref{1.24}) and (\ref{1.30}) should suffice.\\
{\bf Case ab:} We have $z=\sigma_1-\sigma_2-\sigma\le|\sigma_1|+|\sigma_2|+|\sigma|\le 4|\sigma_2|+|\sigma_2|+|\sigma_2|=6|\sigma_2|$ instead of $3|\sigma_2|$ which appears as upper limit of integration and makes no essential difference. We have to check (\ref{1.15}) without use of (\ref{1.7}):  $s-l\le1<2a_1+2a_2-1/2$ by (\ref{1.24})(\ref{1.9}),  $ s-l\le1<2a_2+1/2+2(a+a_1)-2$ by (\ref{1.24})(\ref{1.30}),  $s-l\le1<2a_2+2a-1/2$ by (\ref{1.24})(\ref{1.26}),                    $s-l\le1<2a_2+1/2$ by (\ref{1.24})(\ref{1.25}).\\
{\bf Case ac:} We have $z=\sigma_1-\sigma_2-\sigma\le|\sigma_1|+|\sigma_2|+|\sigma|\le4|\sigma|+|\sigma|+|\sigma|=6|\sigma|$ instead of $3|\sigma|$ which makes no essential difference. In order to check (\ref{1.19}) we only used (\ref{1.24}) instead of (\ref{1.7}).\\
{\bf Region b:} remains unchanged.\\
The proof is complete.
\section{Energy bounds and decomposition of data}
The system (\ref{0.1}),(\ref{0.1a}),(\ref{0.2}) is now transformed into a system of first order in $t$ in the usual way. Defining $ A = -\frac{d^2}{dx^2} $ and
\begin{equation}
\label{2.3a}
n_{\pm} := n \pm iA^{-1/2}n_t
\end{equation}
we have
\begin{eqnarray}
\label{2.3b}
n & = & \frac{1}{2}(n_+ + n_-) \\
\label{2.3c}
2iA^{-1/2}n_t & = & n_+ - n_-
\end{eqnarray}
and the equivalent problem reads as follows
\begin{eqnarray}
\label{2.4}
iu_t + u_{xx}& = &\frac{1}{2}(n_+ + n_-)u \\
in_{\pm t} \mp A^{1/2}n_{\pm}& = &\pm A^{-1/2}(|u|^2)_{xx}
\nonumber
\end{eqnarray}
with initial data
\begin{equation}
\label{2.5}
u(0) = u_0 \quad , \quad n_{\pm}(0) = n_0 \pm iA^{-1/2}n_1
\end{equation}
The standard conservation laws for the original system are: conservation of the $L^2$-norm $\|u(t)\|=: M(u) = M $ and the energy
$$ E:=E(u,n,n_t):= \|u_x(t)\|^2 + 1/2(\|n(t)\|^2 + \|A^{-1/2}n_t(t)\|^2) + \int^{+\infty}_{-\infty} n(t) |u(t)|^2 \, dx $$ 
Now we have by Gagliardo-Nirenberg and $L^2$-conservation:
\begin{eqnarray*}
\left|\int n|u|^2\,dx\right| & \le & \frac{1}{4}\int n^2\,dx + c\int |u|^4\,dx \le \frac{1}{4}\|n\|^2 + c \|u_x\|\,\|u\|^3 \\
& \le & \frac{1}{4}(\|n\|^2 + \|u_x\|^2) + c\|u\|^6 = \frac{1}{4}(\|n\|^2 + \|u_x\|^2) + c\|u_0\|^6 
\end{eqnarray*}
This implies
\begin{equation}
\|u_x(t)\|^2 + \frac{1}{2}(\|n(t)\|^2 + \|A^{-1/2}n_t(t)\|^2) \le E + \frac{1}{4}(\|n(t)\|^2 + \|u_x(t)\|^2) + c_1M^6
\end{equation}
consequently
\begin{eqnarray}
\|u_x(t)\|^2 & \le & \frac{4}{3}(E+c_1M^6)
\label{49a}\\
\|n(t)\|^2 + \|A^{-1/2}n_t(t)\|^2 & \le & 4(E + c_1M^6)
\label{49b}
\end{eqnarray}
We also have
\begin{eqnarray}
E(u,n,n_t) & \le & \|u_x(t)\|^2 + \frac{1}{2}(\|n(t)\|^2 + \|A^{-1/2}n_t(t)\|^2)+ \left| \int^{\infty}_{-\infty} n(t) |u(t)|^2 \,dx \right| \nonumber \\
& \le & \frac{5}{4}\|u_x(t)\|^2 + \frac{3}{4}(\|n(t)\|^2 + \|A^{-1/2}n_t(t)\|^2) + c_1 \|u(t)\|^6
\label{49c}
\end{eqnarray}
These estimates together with $L^2$ - conservation of $u$ and local well-posedness for data in $H^{1,2} \times L^2 \times \dot{H}^{-1,2} $ which is given by \cite{GTV} implies directly also global well-posedness for these data.

Let now data be given with
$$ u_0 \in H^{s,2}({\bf R}) \, , \, n_0 \in L^2({\bf R}) \, , \, n_1 \in \dot{H}^{-1,2}({\bf R}) \quad , \quad 1>s> 9/10 $$
and decompose for $ N \ge 1 $ :
$$ u_0 = u_{01} + u_{02} $$
where
\begin{eqnarray*}
u_{01} & = & {\cal F}^{-1} (\chi_{\{|\xi|\le N\}} \widehat{u_0}(\xi)) = \int_{|\xi|\le N} e^{ix\xi} \widehat{u_0}(\xi)\,d\xi \\
u_{02} & = & {\cal F}^{-1} (\chi_{\{|\xi|\ge N\}} \widehat{u_0}(\xi)) = \int_{|\xi|\ge N} e^{ix\xi} \widehat{u_0}(\xi)\,d\xi
\end{eqnarray*}
One easily shows that
\begin{eqnarray*}
\|u_{01}\|_{H^{l,2}} & \le & c N^{l-s} \|u_0\|_{H^{s,2}} \quad {\mbox for} \, l\ge s \\
\|u_{01}\|_{L^2} & \le & \|u_0\|_{L^2} \\
\|u_{02}\|_{H^{l,2}} & \le & c N^{l-s} \|u_0\|_{H^{s,2}} \quad {\mbox for} \, l\le s \\
\|u_{02}\|_{L^2} & \le & cN^{-s} \|u_0\|_{H^{s,2}}
\end{eqnarray*}
Thus we have the following global bounds for the solution $ (\tilde{u},\tilde{n}) $ of (\ref{0.1}),(\ref{0.1a}) with data $(u_{01},n_0,n_1)$ by (\ref{49c}):
\begin{equation}
\label{50a}
E(\tilde{u},\tilde{n},\tilde{n}_t) \le \frac{5}{4} \|u_{01_x}\|^2 + \frac{3}{4}(\|n_0\|^2 + \|A^{-1/2}n_1\|^2) + c_1\|u_{01}\|^6 \le \bar{c} N^{2(1-s)}
\end{equation}
and thus by (\ref{49a}),(\ref{49b}) and $L^2$-conservation of $\tilde{u}$:
\begin{eqnarray}
\label{2.7}
&& \|\tilde{u}_x(t)\|^2 + \|A^{-1/2}\tilde{n}_t(t)\|^2 + \|\tilde{n}(t)\|^2 \le \hat{c} N^{2(1-s)}  \\
\label{2.8}
&& \|\tilde{u}(t)\| \le M 
\end{eqnarray}
The corresponding global solution $ (\tilde{u},\tilde{n}_{\pm}) $ of (\ref{2.4}) with data $ (u_{01},n_{0+},n_{0-}) $ therefore fulfills:
\begin{eqnarray}
\label{2.9}
&& \|\tilde{u}_x(t)\| \le \hat{c} N^{1-s} \\
\label{2.10}
&& \|\tilde{n}_{\pm}(t)\| \le \hat{c} N^{1-s} \\
&& \|\tilde{u}(t)\| \le M
\label{55}
\end{eqnarray}
where $\hat{c}$ depends essentially only on $\bar{c}$ (the initial energy) and $M$ on the initial $L^2$-norm of $\tilde{u}$.
\section{Further bounds for the regular part}
In order to give further estimates of $ (\tilde{u},\tilde{n}_{\pm}) $ we consider the system of integral equations which belongs to problem (\ref{2.4}) with data $ (u_{01},n_+(0),n_-(0)) $:
\begin{eqnarray}
\label{3.1}
\tilde{u}(t) & = & e^{it\partial_x^2} u_{01} -i \int_0^t e^{i(t-s)\partial_x^2} \frac{1}{2}(\tilde{n}_+(s) + \tilde{n}_-(s))\tilde{u}(s)\,ds \\ \nonumber
\tilde{n}_{\pm}(t) & = & e^{\mp itA^{1/2}} n_{0\pm} \mp i \int_0^t e^{\mp i(t-s)A^{1/2}} A^{-1/2}(|\tilde{u}(s)|^2)_{xx}\,ds 
\label{3.2} 
\end{eqnarray}
We always assume $ t \in I = [0,|I|]$. In this case we could, whenever helpful, place a factor $ \psi_1(t) $ in front of the first terms on the r.h.sides and $ \psi_{|I|}(t) $ in front of any of the integrals in (\ref{3.1}),(\ref{3.2}) without changing the equations at all. Here $\psi \in C_0^{\infty}({\bf R}) $ is a non-negative cut-off function with $ \psi(t) = 0 $ if $ |t|\ge 2 $ , $ \psi(t) = 1 $ if $ |t|\le 1 $ and $\psi_{\delta}:=\psi(t/\delta)$.\\
{\bf Important remark:} Here and in the following section the constants denoted by $c$ or $c_0$ depend essentially only on $\bar{c}$ in (\ref{50a}) (and therefore on $E(\tilde{u},\tilde{n},\tilde{n_t})$ ) and on $M$ (in (\ref{55})).\\
The energy estimate (\ref{2.8}),(\ref{2.9}) gives
\begin{equation}
\label{3.3}
\|\tilde{u}\|_{X^{1,0}(I)} = \|\tilde{u}\|_{L^2(I,H^{1,2}({\bf R}))} \le \|\tilde{u}\|_{L^{\infty}(I,H^{1,2}({\bf R}))} |I|^{1/2} \le c N^{1-s}|I|^{1/2}
\end{equation}
By (\ref{2.10}),(\ref{3.1}),(\ref{0.4}),(\ref{0.5}) and (\ref{0.11}):\\
(in the sequel $a\pm$ denotes a number slightly larger resp. smaller than $a$)
\begin{eqnarray}
\nonumber
\lefteqn{ \|\tilde{u}\|_{X^{0,\frac{1}{2}+}(I)} } \\ \nonumber & \le & c \|u_{01}\|_{L^2({\bf R})} + c(\|\tilde{n}_+\tilde{u}\|_{X^{0,-\frac{1}{2}+}(I)} + \|\tilde{n}_-\tilde{u}\|_{X^{0,-\frac{1}{2}+}(I)})\\
\nonumber
& \le & c + c(\|\tilde{n}_+\tilde{u}\|_{L^{\frac{6}{5}+}(I,L^{\frac{6}{5}+}({\bf R}))} + \|\tilde{n}_-\tilde{u}\|_{L^{\frac{6}{5}+}(I,L^{\frac{6}{5}+}({\bf R}))})\\
\label{3.7}
& \le & c + c\left(\int_I \|\tilde{n}_+\|_{L_x^2}^{\frac{6}{5}+} \|\tilde{u}\|_{L_x^{3+}}^{\frac{6}{5}+}dt + \int_I \|\tilde{n}_-\|_{L_x^2}^{\frac{6}{5}+} \|\tilde{u}\|_{L_x^{3+}}^{\frac{6}{5}+}dt\right)^{(\frac{6}{5}+)^{-1}} \\
\nonumber
& \le & c + c(\|\tilde{n}_+ \|_{L^{\infty}(I,L^2({\bf R}))} + \|\tilde{n}_- \|_{L^{\infty}(I,L^2({\bf R}))})\left(\int_I \|\tilde{u}\|_{L_x^2}^{1-} \|\tilde{u}_x\|_{L_x^2}^{\frac{1}{5}+}dt\right)^{(\frac{6}{5}+)^{-1}} \\
\nonumber
& \le & c + c(\|\tilde{n}_+\|_{L^{\infty}(I,L^2({\bf R}))} + \|\tilde{n}_-\|_{L^{\infty}(I,L^2({\bf R}))}) \|\tilde{u}\|_{L_t^{\infty}(I,L_x^2)}^{\frac{5}{6}-} \|\tilde{u}_x\|_{L_t^{\infty}(I,L_x^2)}^{\frac{1}{6}+} |I|^{\frac{5}{6}-} \\
\nonumber
& \le & c + c(\|\tilde{n}_+\|_{L^{\infty}(I,L^2({\bf R}))} + \|\tilde{n}_-\|_{L^{\infty}(I,L^2({\bf R}))}) N^{\frac{1-s}{6}+} |I|^{\frac{5}{6}-}\\
\nonumber
& \le & c + c N^{1-s} N^{\frac{1-s}{6}+} N^{-\frac{5}{6}4(1-s)+}\\
\nonumber
& \le & 2c
\end{eqnarray}
for $N$ sufficiently large.
We here also used Gagliardo-Nirenberg and (\ref{2.8}),(\ref{2.9}) and assumed $|I| \le cN^{-4(1-s)}$.\\
Next we estimate $ \|\tilde{n}_{\pm}\|_{X^{0,\frac{1}{2}+}(I)} $ in terms of
$ \|\tilde{u}\|_{X^{0,\frac{1}{2}+}(I)} $. From the integral equation (\ref{3.2}),(\ref{0.4}),(\ref{0.5}) we have
\begin{equation}
\|\tilde{n}_{\pm}\|_{X^{0,\frac{1}{2}+}(I)}  \le  c \|n_{0\pm}\|_{L^2({\bf R})} + c \|A^{-1/2}(|\tilde{u}|^2)_{xx}\|_{X^{0,-\frac{3}{8}}(I)} |I|^{\frac{1}{8}-} 
\label{3.8} 
\end{equation} 
By \cite{GTV}, Lemma 4.4 (with $k=1/4,l=0,c=3/8,b_1=3/8$) we have
\begin{equation}
\label{3.9}
\|A^{-1/2}(|\tilde{u}|^2)_{xx}\|_{X^{0,-\frac{3}{8}}(I)} \le c \|\tilde{u}\|^2_{X^{\frac{1}{4},\frac{3}{8}}(I)} \le c \|\tilde{u}\|_{X^{0,\frac{1}{2}+}(I)}^{\frac{3}{2}} \|\tilde{u}\|_{X^{1,0}(I)}^{\frac{1}{2}} 
\end{equation}
where the last estimate follows by interpolation from (\ref{0.9}). 
From (\ref{3.3}),(\ref{3.7}),(\ref{3.8}),(\ref{3.9}) we get
\begin{equation}
\|\tilde{n}_{\pm}\|_{X^{0,\frac{1}{2}+}(I)} \le c \|n_{0\pm}\|_{L^2({\bf R})}+ c\|\tilde{u}\|_{X^{0,\frac{1}{2}+}(I)}^{\frac{3}{2}} N^{\frac{1-s}{2}} |I|^{\frac{3}{8}-} \le c(N^{1-s} + N^{\frac{1-s}{2}})
\label{3.10}
\end{equation}
Thus we have
\begin{lemma}
\label{Lemma 3.1}
Let $ |I| \le N^{-4(1-s)} $ and $\|n_{0\pm}\|_{L^2({\bf R})} \le cN^{1-s} $. Then we have
\begin{eqnarray}
\label{3.11}
\|\tilde{n}_{\pm}\|_{X^{0,\frac{1}{2}+}(I)} & \le & cN^{1-s} \\
\|\tilde{u}\|_{X^{0,\frac{1}{2}+}(I)} & \le & c
\label{3.12}
\end{eqnarray}
\end{lemma}
The next step is an estimate of $ \|\tilde{u}\|_{X^{1,\frac{1}{2}}(I)} $.\\ From the integral equation (\ref{3.1}) we get by ({\ref{0.4}),(\ref{0.6}):
\begin{eqnarray} \nonumber
\|\tilde{u}\|_{X^{1,\frac{1}{2}}(I)} & \le & c \|u_{01}\|_{H^{1,2}({\bf R})} + c (\|\tilde{n}_+\tilde{u}\|_{X^{1,-\frac{1}{2}}(I)} + \|\tilde{n}_-\tilde{u}\|_{X^{1,-\frac{1}{2}}(I)}) \\ \label{3.14} && + c(\|\tilde{n}_+\tilde{u}\|_{Y^1(I)} + \|\tilde{n}_-\tilde{u}\|_{Y^1(I)}) \\ \nonumber
& \le & c N^{1-s} + c(\|\tilde{n}_+\|_{X^{0,\frac{1}{2}+}(I)} + \|\tilde{n}_-\|_{X^{0,\frac{1}{2}+}(I)}) \|\tilde{u}\|_{X^{\frac{1}{2}+,\frac{1}{4}+}(I)} \\ \nonumber
& \le & c N^{1-s} + c(\|\tilde{n}_+\|_{X^{0,\frac{1}{2}+}(I)} + \|\tilde{n}_-\|_{X^{0,\frac{1}{2}+}(I)}) \|\tilde{u}\|_{X^{0,\frac{1}{2}+}(I)}^{\frac{1}{2}-} \|\tilde{u}\|_{X^{1,0}(I)}^{\frac{1}{2}+} \\ 
& \le & c N^{1-s} + c(\|\tilde{n}_+\|_{X^{0,\frac{1}{2}+}(I)} + \|\tilde{n}_-\|_{X^{0,\frac{1}{2}+}(I)}) \|\tilde{u}\|_{X^{0,\frac{1}{2}+}(I)}^{\frac{1}{2}-} N^{\frac{1-s}{2}+}|I|^{\frac{1}{4}+}  
\nonumber    
\end{eqnarray}
Here we used Prop. \ref{Proposition 1.1} and Prop. \ref{Proposition1.2} (with $s=1,l=0,k=\frac{1}{2}+,a_1=\frac{1}{2},a=\frac{1}{2}+,a_2=\frac{1}{4}+$), an interpolation argument and (\ref{3.3}). Consequently we get
\begin{lemma}
\label{Lemma 3.3}
If $|I| \le N^{-4(1-s)} $ and $ \|n_{0\pm}\|_{L^2({\bf R})} \le c N^{1-s} $, the following estimate holds
\begin{equation}
\label{3.16}
\|\tilde{u}\|_{X^{1,\frac{1}{2}}(I)} \le cN^{1-s} 
\end{equation}
\end{lemma}
{\bf Proof:} follows immediately from Lemma \ref{Lemma 3.1} and (\ref{3.14}):
$$\|\tilde{u}\|_{X^{1,\frac{1}{2}}(I)} \le cN^{1-s} + cN^{1-s}N^{\frac{1-s}{2}+}N^{-4(1-s)\frac{1}{4}+} \le cN^{1-s} $$
{\bf Remark:} If the data fulfill the conditions
\begin{eqnarray}
\label{3.18a}
\|u_{01}\|_{L^2({\bf R})} & \le & c \\
\label{3.18b}
\|u_{01_x}\|_{L^2({\bf R})} & \le &  cN^{1-s} \\
\label{3.18c}
\|n_{0+}\|_{L^2({\bf R})} + \|n_{0-}\|_{L^2({\bf R})}& \le & cN^{1-s} 
\end{eqnarray}
then the following estimates hold on $|I|\le N^{-4(1-s)}$:
\begin{eqnarray}
&& \|\tilde{n}_{\pm}\|_{X^{0,\frac{1}{2}+}(I)} \le cN^{1-s}
\label{76} \\
&& \|\tilde{u}\|_{X^{1,\frac{1}{2}}(I)} \le cN^{1-s} \\
&&\|\tilde{u}\|_{X^{0,\frac{1}{2}+}(I)} \le c
\end{eqnarray}
This follows from Lemma \ref{Lemma 3.1}, \ref{Lemma 3.3}.\\
Also the estimates (\ref{2.9}),(\ref{2.10}),(\ref{55}) hold under these assumptions.
\section{The part with rough data}
Let $(u,n_+,n_-)$ be a solution of (\ref{2.4}) with data $(u_0,n_{0+},n_{0-})$ and $(\tilde{u},\tilde{n}_+,\tilde{n}_-)$ the solution with data $(u_{01},n_{0+},n_{0-})$.\\
Define $ v:= u-\tilde{u} \, , \, m_{\pm} := n_{\pm} - \tilde{n}_{\pm}\,.$  Then $ (v,m_+,m_-) $ fulfills
\begin{eqnarray}
\nonumber
iv_t + v_{xx} & = & iu_t+u_{xx}-i\tilde{u}_t-\tilde{u}_{xx} = \frac{1}{2}(n_++n_-)u-\frac{1}{2}(\tilde{n}_++\tilde{n}_-)\tilde{u} \\
\nonumber
& = & \frac{1}{2}(\tilde{n}_++m_++\tilde{n}_-+m_-)(\tilde{u}+v) - \frac{1}{2}(\tilde{n}_++\tilde{n}_-)\tilde{u} \\ \nonumber
& = & \frac{1}{2}(\tilde{n}_++\tilde{n}_-)v + \frac{1}{2}(m_++m_-)v + \frac{1}{2}(m_++m_-)\tilde{u} \\ \label{4.1}
& = & F_1 + F_2 + F_3 =: F
\end{eqnarray}
and
\begin{eqnarray}
\nonumber
&& im_{\pm t} \mp A^{1/2}m_{\pm}  =  in_{\pm t}-i\tilde{n}_{\pm t} \mp A^{1/2}n_{\pm} \pm A^{1/2}\tilde{n}_{\pm} \\
\nonumber
& & = \pm A^{-1/2}(|u|^2)_{xx} \mp A^{-1/2}(|\tilde{u}|^2)_{xx}  \\
\nonumber
& &= \pm A^{-1/2}((\tilde{u}+v)(\bar{\tilde{u}}+\bar{v}))_{xx} \mp A^{-1/2}(\tilde{u}\bar{\tilde{u}})_{xx}  \\
 \nonumber
&  &= \pm A^{-1/2}(\tilde{u}\bar{v})_{xx} \pm A^{-1/2}(|v|^2)_{xx} \pm A^{-1/2}(v\bar{\tilde{u}})_{xx}  \\
\nonumber 
&  &=: G_1 + G_2 + G_3  \\ \label{4.2} && =:  G 
\end{eqnarray}
Furthermore
\begin{eqnarray}
\label{4.3}
&& v(0) = u(0)-\tilde{u}(0) = u_0 - u_{01} = u_{02} \\
&& m_{\pm}(0) = n_{\pm}(0) - \tilde{n}_{\pm}(0) = 0
\end{eqnarray}
The corresponding system of integral equations reads as follows
\begin{eqnarray}
\label{4.5}
v(t) & = & e^{it\partial_x^2} u_{02} -i \int_0^t e^{i(t-s)\partial_x^2}F(s)\,ds \\
m_{\pm}(t) & = & -i \int_0^t e^{\mp i(t-s)A^{1/2}} G(s) ds 
\label{4.6}
\end{eqnarray}
Here we have $ u_{02} \in H^{s,2}({\bf R}) $ with
\begin{eqnarray}
&& \|u_{02}\|_{H^{s,2}} \le c \|u_0\|_{H^{s,2}} \le c 
\label{4.7} \\
\label{4.8}
&& \|u_{02}\|_{L^2} \le cN^{-s} \|u_0\|_{H^{s,2}} \le cN^{-s}
\end{eqnarray}
We construct a solution of (\ref{4.5}),(\ref{4.6}) in some time interval $I$ by the standard contraction mapping principle. We define
\begin{equation}
\label{4.8a}
w(t):=-i \int_0^t e^{i(t-s)\partial_x^2}F(s)\,ds \quad , \quad z_{\pm}(t):= \mbox{r.h.s. of (\ref{4.6})}
\end{equation}
and a mapping $S=(S_0,S_+,S_-)$ by
\begin{eqnarray*}
(S_0v)(t) & := & e^{it\partial_x^2} u_{02} + w(t) \\
(S_{\pm}m_{\pm})(t) & := & z_{\pm}(t)
\end{eqnarray*}
\begin{prop}
\label{Proposition 4.1}
For $9/10 < s < 1$ and given data $u_{02} \in H^{s,2}({\bf R}) $ with (\ref{4.7}),(\ref{4.8}) and $u_{01},n_{0\pm}$ as in (\ref{3.18a}),(\ref{3.18b}),(\ref{3.18c}) the system of integral equations (\ref{4.5}),(\ref{4.6}) has a unique solution $(v,m_{\pm})\in X^{s,\frac{1}{2}+}(I) \times X^{0,\frac{1}{2}+}(I)$ in the same interval $I$ with $|I|=N^{-4(1-s)-\delta}$, $\delta > 0$ of the preceding section, which fulfills
\begin{eqnarray}
\label{4.9}
\|v\|_{X^{0,\frac{1}{2}+}(I)} & \le & cN^{-s} \\
\label{4.10}
\|v\|_{X^{s,\frac{1}{2}+}(I)} & \le & c \\
\label{4.11}
\|m_{\pm}\|_{X^{-\frac{1}{2},\frac{1}{2}+}(I)} & \le & cN^{-s} \\
\label{4.12}
\|m_{\pm}\|_{X^{0,\frac{1}{2}+}(I)} & \le & cN^{-\frac{1}{2}-\frac{1}{4}s-\frac{\delta}{4}+}
\end{eqnarray}
\end{prop}
{\bf Proof:} We want to use Banach's fixed point theorem in the set $Z$, where
\begin{eqnarray*}
Z&:=&\{ \|v\|_{X^{0,\frac{1}{2}+\epsilon}(I)} \le c_0N^{-s}\, , \,\|v\|_{X^{s,\frac{1}{2}+\epsilon}(I)} \le c_0 \\ && \;\,\|m_{\pm}\|_{X^{-\frac{1}{2},\frac{1}{2}+\epsilon}(I)} \le c_0 N^{-s} \,, \,\|m_{\pm}\|_{X^{0,\frac{1}{2}+\epsilon}(I)} \le c_0 N^{-\frac{1}{2}-\frac{1}{4}s-\frac{\delta}{4}+} \} 
\end{eqnarray*} 
with its natural metric, $c_0$ chosen below.\\
Now take any $(v,m_+,m_-)\in Z$. In order to show $(S_0v,S_+m_+,S_-m_-)\in Z$ we estimate $\|S_0v\|_{X^{s,\frac{1}{2}+\epsilon}(I)}$ first.\\
We have by Prop. \ref{Proposition 1.1} (with $k=s-\frac{1}{2}+8\epsilon,l=0,a_1=\frac{1}{2}-2\epsilon,a=a_2=\frac{1}{2}+\epsilon$) and interpolation
\begin{eqnarray*}
\|\tilde{n}_{\pm}v\|_{X^{s,-\frac{1}{2}+2\epsilon}(I)} & \le & c \|\tilde{n}_{\pm}\|_{X^{0,\frac{1}{2}+\epsilon}(I)} \|v\|_{X^{s-\frac{1}{2}+8\epsilon,\frac{1}{2}+\epsilon}(I)} \\
& \le & c\|\tilde{n}_{\pm}\|_{X^{0,\frac{1}{2}+\epsilon}(I)} \|v\|_{X^{0,\frac{1}{2}+\epsilon}(I)}^{\frac{1}{2s}-} \|v\|_{X^{s,\frac{1}{2}+\epsilon}(I)}^{1-\frac{1}{2s}+} \\
& \le & cN^{1-s}N^{-\frac{1}{2}+} \le cN^{\frac{1}{2}-s+} \\
\|m_{\pm}v\|_{X^{s,-\frac{1}{2}+2\epsilon}(I)} & \le & c \|m_{\pm}\|_{X^{0,\frac{1}{2}+\epsilon}(I)} \|v\|_{X^{s-\frac{1}{2}+8\epsilon,\frac{1}{2}+\epsilon}(I)} \\
& \le & c\|m_{\pm}\|_{X^{0,\frac{1}{2}+\epsilon}(I)} \|v\|_{X^{0,\frac{1}{2}+\epsilon}(I)}^{\frac{1}{2s}-} \|v\|_{X^{s,\frac{1}{2}+\epsilon}(I)}^{1-\frac{1}{2s}+} \\
& \le & cN^{-\frac{1}{2}-\frac{1}{4}s+}N^{-\frac{1}{2}+} = cN^{-1-\frac{1}{4}s+} \\
\|m_{\pm}\tilde{u}\|_{X^{s,-\frac{1}{2}+2\epsilon}(I)} & \le & c \|m_{\pm}\|_{X^{0,\frac{1}{2}+\epsilon}(I)} \|\tilde{u}\|_{X^{s-\frac{1}{2}+8\epsilon,\frac{1}{2}+\epsilon}(I)} \\
& \le & c\|m_{\pm}\|_{X^{0,\frac{1}{2}+\epsilon}(I)} \|\tilde{u}\|_{X^{0,\frac{1}{2}+\epsilon}(I)}^{\frac{3}{2}-s-} \|\tilde{u}\|_{X^{1-,\frac{1}{2}+\epsilon}(I)}^{s-\frac{1}{2}+} \\
& \le & cN^{-\frac{1}{2}-\frac{1}{4}s}N^{(1-s)(s-\frac{1}{2})+} = cN^{-s(s-\frac{1}{2})-1+\frac{3}{4}s+}
\end{eqnarray*}
We have $ 1>s> 9/10 $, so the exponents of $N$ are negative. Thus with $\gamma(s)>0$:
$$ \|F\|_{X^{s,-\frac{1}{2}+2\epsilon}(I)} \le c N^{-\gamma(s)} $$ and therefore with $c_0 \ge 2c\|u_{02}\|_{H^{s,2}({\bf R})} $ we have:
$$ \|S_0v\|_{X^{s,\frac{1}{2}+\epsilon}(I)} \le c\|u_{02}\|_{H^{s,2}({\bf R})} + c\|F\|_{X^{s,-\frac{1}{2}+2\epsilon}(I)} \le \frac{c_0}{2} + cN^{-\gamma(s)} \le c_0 $$
if $N$ is sufficiently large.\\
Next we estimate $ \|S_0v\|_{X^{0,\frac{1}{2}+\epsilon}(I)} $. We use \cite{GTV},Lemma 4.3 (with $k=0$ , $l=-\frac{1}{2}$ , $c_1=\frac{1}{2}-$ , $b=\frac{1}{2}+$ , $b_1=\frac{1}{2}+$ and with $k=l=0$ , $c_1=\frac{1}{4}+$ , $b=b_1=\frac{1}{2}+$ ) and get:
\begin{eqnarray*}
\lefteqn{\|S_0v\|_{X^{0,\frac{1}{2}+\epsilon}(I)} }\\
& \le & c\|u_{02}\|_{L^2({\bf R})} \\
&& +c ( \|m_{\pm}v\|_{X^{0,-\frac{1}{2}+2\epsilon}(I)}|I|^{\epsilon} + \|\tilde{n}_{\pm}v\|_{X^{0,-\frac{1}{4}-}(I)}|I|^{\frac{1}{4}-}
+ \|m_{\pm}\tilde{u}\|_{X^{0,-\frac{1}{2}+2\epsilon}(I)}|I|^{\epsilon}) \\
& \le & c_1 N^{-s} +c (\|m_{\pm}\|_{X^{-\frac{1}{2},\frac{1}{2}+}(I)}\|v\|_{X^{0,\frac{1}{2}+}(I)}|I|^{\epsilon}
+ \|\tilde{n}_{\pm}\|_{X^{0,\frac{1}{2}+}(I)}\|v\|_{X^{0,\frac{1}{2}+}(I)}|I|^{\frac{1}{4}-} \\ 
&& +  \|m_{\pm}\|_{X^{-\frac{1}{2},\frac{1}{2}+}(I)}\|\tilde{u}\|_{X^{0,\frac{1}{2}+}(I)}|I|^{\epsilon}) \\ 
& \le & c_1 N^{-s} + c(N^{-s}N^{-s}|I|^{\epsilon} + N^{1-s} N^{-s}|I|^{\frac{1}{4}-} + N^{-s}|I|^{\epsilon} ) \\
& \le & c(N^{-s}|I|^{\epsilon} + N^{1-2s} |I|^{\frac{1}{4}-}) + c_1 N^{-s} \\
& \le & c_1N^{-s} + cN^{-s}|I|^{\epsilon} + cN^{1-2s-(1-s)-\frac{\delta}{4}+} \\
& \le & N^{-s}(\frac{c_0}{2} + c N^{-4(1-s)\epsilon} + cN^{-\frac{\delta}{4}+}) \\
& \le & c_0 N^{-s}
\end{eqnarray*}
from the definition of $Sv$,(\ref{4.8}),(\ref{0.4}) and (\ref{0.5}), where $ c_0 \ge 2c_1 $ , $\epsilon>0$ small and $N$ sufficiently large.\\
Next we treat $ \|S_{\pm}m_{\pm}\|_{X^{0,\frac{1}{2}+\epsilon}(I)} $. By \cite{GTV}, Lemma 4.4 (with $l=0$ , $ k=\frac{1}{4}$ , $c=\frac{1}{4}+$ , $b_1=\frac{1}{2}$):
\begin{eqnarray*}
\lefteqn{ \|(\tilde{u}\bar{v})_x\|_{X^{0,-\frac{1}{4}-}(I)} + \|(|v|^2)_x\|_{X^{0,-\frac{1}{4}-}(I)} + \|(v\bar{\tilde{u}})_x\|_{X^{0,-\frac{1}{4}-}(I)} } \\
&& \le c (\|\tilde{u}\|_{X^{\frac{1}{4},\frac{1}{2}}(I)}\|v\|_{X^{\frac{1}{4},\frac{1}{2}}(I)} + \|v\|_{X^{\frac{1}{4},\frac{1}{2}}(I)}^2)\\
&& \le c (\|\tilde{u}\|_{X^{0,\frac{1}{2}}(I)}^{\frac{3}{4}-}\|\tilde{u}\|_{X^{1-,\frac{1}{2}}(I)}^{\frac{1}{4}+}\|v\|_{X^{0,\frac{1}{2}+}(I)}^{1-\frac{1}{4s}}\|v\|_{X^{s,\frac{1}{2}+}(I)}^{\frac{1}{4s}} + \|v\|_{X^{0,\frac{1}{2}+}(I)}^{2(1-\frac{1}{4s})}\|v\|_{X^{s,\frac{1}{2}+}(I)}^{\frac{1}{2s}}) \\
&& \le c(N^{\frac{1-s}{4}+}N^{-s(1-\frac{1}{4s})}+ N^{-2s(1-\frac{1}{4s})}) \\
&& = c(N^{\frac{1}{2}-\frac{5}{4}s+} + N^{\frac{1}{2}-2s}) \\
&& \le cN^{\frac{1}{2}-\frac{5}{4}s+}
\end{eqnarray*}
Thus
\begin{eqnarray}
&& \|S_{\pm}m_{\pm}\|_{X^{0,\frac{1}{2}+}(I)} = \| \int_0^t e^{\mp i(t-s)A^{1/2}} G(s)\,ds\|_{X^{0,\frac{1}{2}+}(I)}
 \le c\|G\|_{X^{0,-\frac{1}{4}-}(I)}|I|^{\frac{1}{4}-}
\nonumber \\
&& \le cN^{\frac{1}{2}-\frac{5}{4}s+}|I|^{\frac{1}{4}-} = c N^{\frac{1}{2}-\frac{5}{4}s-(1-s)-\frac{\delta}{4}+} = c N^{-\frac{1}{2}-\frac{1}{4}s-\frac{\delta}{4}+} 
 \le c_0 N^{-\frac{1}{2}-\frac{1}{4}s-\frac{\delta}{4}+}
\label{93} 
\end{eqnarray}
if $N$ is chosen sufficiently large.\\
Finally we treat $\|S_{\pm}m_{\pm}\|_{X^{-\frac{1}{2},\frac{1}{2}+\epsilon}(I)}$. By \cite{GTV}, Lemma 4.4 (with $l=-\frac{1}{2}$ , $ k=0$ , $b_1=\frac{1}{2}+$ , $ c=\frac{1}{4}+$):
\begin{eqnarray*}
\|G\|_{X^{-\frac{1}{2},-\frac{1}{4}-}(I)} & \le & \|(\tilde{u}\bar{v})_x\|_{X^{-\frac{1}{2},-\frac{1}{4}-}(I)} + \|(|v|^2)_x\|_{X^{-\frac{1}{2},-\frac{1}{4}-}(I)} + \|(v\bar{{\tilde{u}}})_x\|_{X^{-\frac{1}{2},-\frac{1}{4}-}(I)} \\
& \le & c(\|\tilde{u}\|_{X^{0,\frac{1}{2}+}(I)}
 \|v\|_{X^{0,\frac{1}{2}+}(I)} + \|v\|_{X^{0,\frac{1}{2}+}(I)}^2) \\
& \le & c(N^{-s} + N^{-2s}) \\
& \le & cN^{-s}
\end{eqnarray*}
Thus
\begin{eqnarray}
\lefteqn{\|S_{\pm}m_{\pm}\|_{X^{-\frac{1}{2},\frac{1}{2}+\epsilon}(I)} = \| \int_0^t e^{\mp i(t-s)A^{1/2}}G(s)ds\|_{X^{-\frac{1}{2},\frac{1}{2}+\epsilon}(I)} }
\nonumber \\
&& \le c\|G\|_{X^{-\frac{1}{2},-\frac{1}{4}-}(I)} |I|^{\frac{1}{4}-} \le cN^{-s}|I|^{\frac{1}{4}-} \le c_0 N^{-s}
\label{94}
\end{eqnarray}  
for $N$ sufficiently large.\\
Summarizing we have shown that $S$ maps $Z$ into itself.\\
The contraction property uses exactly the same type of estimates and its proof is therefore omitted. Thus the proposition is proved.

The next estimates show that the nonlinear part $w$ of (\ref{4.5}), defined by (\ref{4.8a}), behaves better than the linear part.\\
We use (\ref{4.14b}) in connection with Prop. \ref{Proposition1.2} (with $s=1,l=0,k=\frac{1}{2}+,a=a_2=\frac{1}{2}$). This gives
\begin{eqnarray}
\lefteqn{ \|w\|_{L^{\infty}(I,H^{1,2}({\bf R}))} } \nonumber \\ \nonumber
&& \le c\left(\|(m_++m_-)v\|_{Y^1(I)} + \|(\tilde{n}_++\tilde{n}_-)v\|_{Y^1(I)}
+ \|(m_++m_-)\tilde{u}\|_{Y^1(I)}\right) \\ \nonumber
&& \le c(\|m_++m_-\|_{X^{0,\frac{1}{2}}(I)}\|v\|_{X^{\frac{1}{2}+,\frac{1}{2}}(I)} + \|\tilde{n}_++\tilde{n}_-\|_{X^{0,\frac{1}{2}}(I)}\|v\|_{X^{\frac{1}{2}+,\frac{1}{2}}(I)} \\ \nonumber
&& \qquad + \|m_++m_-\|_{X^{0,\frac{1}{2}}(I)}\|\tilde{u}\|_{X^{\frac{1}{2}+,\frac{1}{2}}(I)})\\ \nonumber
&& \le c(\|m_++m_-\|_{X^{0,\frac{1}{2}}(I)}\|v\|_{X^{0,\frac{1}{2}}(I)}^{1-\frac{1}{2s}-}\|v\|_{X^{s,\frac{1}{2}}(I)}^{\frac{1}{2s}+} \\ \nonumber
&& \qquad + \|\tilde{n}_++\tilde{n}_-\|_{X^{0,\frac{1}{2}}(I)}\|v\|_{X^{0,\frac{1}{2}}(I)}^{1-\frac{1}{2s}-}\|v\|_{X^{s,\frac{1}{2}}(I)}^{\frac{1}{2s}+} \\ \nonumber
&& \qquad + \|m_++m_-\|_{X^{0,\frac{1}{2}}(I)}\|\tilde{u}\|_{X^{0,\frac{1}{2}}(I)}^{\frac{1}{2}-}\|\tilde{u}\|_{X^{1-,\frac{1}{2}}(I)}^{\frac{1}{2}+}) \\ \nonumber
&& \le c(N^{-\frac{1}{2}-\frac{1}{4}s} N^{-s(1-\frac{1}{2s})+} + N^{1-s}N^{-s(1-\frac{1}{2s})+} +  N^{-\frac{1}{2}-\frac{1}{4}s+}N^{\frac{1-s}{2}+}) \\ \label{95}
&& = c(N^{-\frac{5}{4}s} + N^{\frac{3}{2}-2s+} + N^{-\frac{3}{4}s} ) \le c N^{\frac{3}{2}-2s+}
\end{eqnarray}
because $ \frac{3}{2}-2s > -\frac{3}{4}s $ , 
by (\ref{93}),(\ref{4.9}),(\ref{4.10}),(\ref{3.11}),(\ref{3.12}),(\ref{3.16}), using $|I| = N^{-4(1-s)-}$. 
The next step is to estimate $\|w\|_{L^{\infty}(I,L^2({\bf R}))}$.\\
We use \cite{GTV}, Lemma 4.3 (with $k=l=0$ , $ b=b_1=\frac{1}{2}+$ , $c_1=\frac{1}{4}+$):
\begin{eqnarray}
\nonumber
\lefteqn{ \|w\|_{L^{\infty}(I,L^2({\bf R}))} \le c \|w\|_{X^{0,\frac{1}{2}+}(I)} } \\
\nonumber && \le c(\|(m_++m_-)v\|_{X^{0,-\frac{1}{4}-}(I)} + \|(\tilde{n}_++\tilde{n}_-)v\|_{X^{0,-\frac{1}{4}-}(I)} \\ 
\nonumber && \qquad
+ \|(m_++m_-)\tilde{u}\|_{X^{0,-\frac{1}{4}-}(I)})|I|^{\frac{1}{4}-} \\
\nonumber && \le c(\|m_++m_-\|_{X^{0,\frac{1}{2}+}(I)}\|v\|_{X^{0,\frac{1}{2}+}(I)} + \|\tilde{n}_++\tilde{n}_-\|_{X^{0,\frac{1}{2}+}(I)}\|v\|_{X^{0,\frac{1}{2}+}(I)} \\
\nonumber && \qquad
+ \|m_++m_-\|_{X^{0,\frac{1}{2}+}(I)}\|\tilde{u}\|_{X^{0,\frac{1}{2}+}(I)})|I|^{\frac{1}{4}-} \\
\nonumber
&& \le c(N^{-\frac{1}{2}-\frac{1}{4}s}N^{-s} + N^{1-s} N^{-s} + N^{-\frac{1}{2}-\frac{1}{4}s})|I|^{\frac{1}{4}-} \\
\nonumber && \le c(N^{1-2s} + N^{-\frac{1}{2}-\frac{1}{4}s})N^{-(1-s)-} \\ 
\label{4.16}
&& \le c N^{-\frac{3}{2}+\frac{3}{4}s+} 
\end{eqnarray}
by (\ref{94}),(\ref{4.9}),(\ref{76}),(\ref{3.12}), because $ 1-2s < -\frac{1}{2}-\frac{1}{4}s \Leftrightarrow s > \frac{6}{7}$\\
We also repeat the corresponding estimate for $m_{\pm} = Sm_{\pm} = z_{\pm} $ from (\ref{93}):
\begin{equation}
\label{4.17} 
\|z_{\pm}\|_{X^{0,\frac{1}{2}+}(I)} = \|m_{\pm}\|_{X^{0,\frac{1}{2}+}(I)}  \le
 cN^{-\frac{1}{2}-\frac{1}{4}s-} 
\end{equation}
\section{The iteration process}
In the preceding sections we have constructed a solution $(u,n)$ of our original problem (\ref{0.1}),(\ref{0.1a}) with data (\ref{0.2}) $(u_0,n_0,n_1)$ in the time interval $I=[0,|I|]$ with $|I|=N^{-4(1-s)-}$. Namely, if we define $u:=v+\tilde{u}$ , $ n_{\pm}:=m_{\pm}+\tilde{n}_{\pm}$ we easily see by (\ref{4.1}),(\ref{4.2}) that $(u,n_+,n_-)$ satisfies the system (\ref{2.4}). Moreover the initial conditions $u(0)=u_0$ , $ n_{\pm}(0)=n_{0\pm}$ are satisfied. This initial value problem is equivalent to the original system (\ref{0.1}),(\ref{0.1a}) by (\ref{2.3b}). The initial data are transformed via (\ref{2.3b}),(\ref{2.3c}) by $n_{0\pm}=\frac{1}{2}(n_0+n_1)$ , $ 2iA^{-1/2}n_1=n_{0+}-n_{0-}$ or conversely by $n_{0\pm}=n_0 \pm iA^{-1/2}n_1$.\\
In order to continue the solution of (\ref{2.4}),(\ref{2.5}) we take as new initial data for our system (\ref{2.4}) the triple $(\tilde{u}(|I|)+w(|I|),\tilde{n}_+(|I|)+z_+(|I|),\tilde{n}_-(|I|)+z_-(|I|))$ instead of $(u_{01},n_{0+},n_{0-})$. When we have shown that this problem has a solution $(\tilde{\tilde{u}},\tilde{\tilde{n}}_+,\tilde{\tilde{n}}_-)$ in the time interval $[|I|,2|I|]$ of equal length $|I|$ we insert this solution into the system (\ref{4.1}),(\ref{4.2}) in place of $(\tilde{u},\tilde{n}_+,\tilde{n}_-)$ and solve this problem with data $(e^{i|I|\partial_x^2}u_{02},0,0)$ in $[|I|,2|I|]$. Adding up the solutions we get a solution of the original problem in $[|I|,2|I|]$ as before. This defines an iteration process. At each step we have to ensure the same bounds on the initial data which were used in the first step. The replacement of $u_{02}$ by$e^{i|I|\partial_x^2}u_{02}$ is harmless, because $ e^{i|I|\partial^2_x} $ is unitary in $\dot{H}^s({\bf R})$. These bounds are controlled by the energy and the $L^2$-conservation law (cf. (\ref{49a}),(\ref{49b})). Thus we have to estimate these quantities independently of the iteration step. This is easy for $L^2$-conservation, the increment when replacing $u_{01}$ by $u_{01}(|I|)+w(|I|)$ using $L^2$-conservation is given by
\begin{eqnarray*}
&& \left| \|\tilde{u}(|I|)+w(|I|)\|_{L^2({\bf R})} - \|u_{01}\|_{L^2({\bf R})} \right| =
 \left| \|\tilde{u}(|I|)+w(|I|)\|_{L^2({\bf R})} - \|\tilde{u}(|I|)\|_{L^2({\bf R})} \right| \\
&& \le \|w(|I|)\|_{L^2({\bf R})} \le c_2 N^{-\frac{3}{2}+\frac{3}{4}s+}
\end{eqnarray*}
by (\ref{4.16}), where $c_2 = c_2(\bar{c},M) $.\\
The number of iteration steps in order to reach the given time $T$ is $ \frac{T}{|I|} = TN^{4(1-s)+} $. This means that we have to ensure in order to get uniform control over the $L^2$-norm of $\tilde{u},\tilde{\tilde{u}}, \dots $:
$$ c_2TN^{4(1-s)+}+ N^{-\frac{3}{2}+\frac{3}{4}s+} < M $$
where $ c_2 = c_2(2\bar{c},2M) $ (remark that initially the $L^2$-norm of $\tilde{u}$ was also bounded by $M$).\\
This is fulfilled for $N$ sufficiently large if $ 4(1-s)-\frac{3}{2}+\frac{3}{4}s < 0 \Leftrightarrow s > \frac{10}{13} $ which is fulfilled.\\
The increment of the energy is given by
\begin{eqnarray*}
&& \hspace{-0.7cm}  \left| E(\tilde{u}(|I|)+w(|I|),n(|I|)+m(|I|),n_t(|I|)+m_t(|I|)) - E(u_{01},n_0,n_1) \right| \\
&& \hspace{-0.7cm}= \left| E(\tilde{u}(|I|)+w(|I|),n(|I|)+m(|I|),n_t(|I|)+m_t(|I|)) - E(\tilde{u}(|I|),n(|I|),n_t(|I|)) \right| \\
&& \hspace{-0.7cm}\le 2 (\|\tilde{u}_x(|I|)\| + \|w_x(|I|)\|)\|w_x(|I|)\| + (\|n(|I|)\| + \|m(|I|)\|)\|w(|I|)\| \\
&& \hspace{-0.7cm}+ (\|A^{-1/2}n_t(|I|)\| + \|A^{-1/2}m_t(|I|)\|)\|A^{-1/2}m_t(|I|)\| \\
&& \hspace{-0.7cm}+ \int^{\infty}_{-\infty} |m(|I|)||\tilde{u}(|I|)+w(|I|)|^2 \, dx + \int^{\infty}_{-\infty} |n(|I|)|||\tilde{u}(|I|)+w(|I|)|^2 - |\tilde{u}(|I|)|^2| \, dx
\end{eqnarray*}
Using (\ref{2.7}),(\ref{95}) the first term is bounded by
$$ c(N^{1-s} + N^{\frac{3}{2}-2s+})N^{\frac{3}{2}-2s+} \le cN^{1-s} N^{\frac{3}{2}-2s+} $$
the second and third one using (\ref{2.7}),(\ref{4.17}) by
$$ c(N^{1-s}+N^{-\frac{1}{2}-\frac{1}{4}s-})N^{-\frac{1}{2}-\frac{1}{4}s-} \le cN^{1-s}N^{-\frac{1}{2}-\frac{1}{4}s-} $$
The fourth term is estimated using Gagliardo-Nirenberg and (\ref{2.7}),(\ref{95}),(\ref{4.16}),(\ref{4.17}):
\begin{eqnarray*}
&& \int^{\infty}_{-\infty} |m(|I|)||\tilde{u}(|I|)+w(|I|)|^2 \, dx \le 2 \int^{\infty}_{-\infty} |m(|I|)|(|\tilde{u}(|I|)|^2+|w(|I|)|^2) \, dx \\
&& \le 2\|m(|I|)\|_{L^2({\bf R})}( \|\tilde{u}(|I|)\|_{L^2({\bf R})} \|\tilde{u}(|I|)\|_{L^{\infty}({\bf R})}+ \|w(|I|)\|_{L^2({\bf R})} \|w(|I|)\|_{L^{\infty}({\bf R})}) \\
&& \le c\|m(|I|)\|_{L^2({\bf R})}(\|\tilde{u}(|I|)\|_{L^2({\bf R})}^{\frac{3}{2}} \|\tilde{u}_x(|I|)\|_{L^2({\bf R})}^{\frac{1}{2}} + \|w(|I|)\|_{L^2({\bf R})}^{\frac{3}{2}} \|w_x(|I|)\|_{L^2({\bf R})}^{\frac{1}{2}}) \\
&& \le cN^{-\frac{1}{2}-\frac{1}{4}s-} (N^{\frac{1-s}{2}} + N^{\frac{3}{2}(-\frac{3}{2}+\frac{3}{4}s)+} N^{\frac{1}{2}(\frac{3}{2}-2s)+}) \\
&& \le cN^{-\frac{1}{2}-\frac{1}{4}s-}N^{\frac{1-s}{2}} = cN^{-\frac{3}{4}s-}
\end{eqnarray*}
The fifth term is similarly estimated as follows:
\begin{eqnarray*}
\lefteqn{\int^{\infty}_{-\infty} |n(|I|)|||\tilde{u}(|I|)+w(|I|)|^2 - |\tilde{u}(|I|)|^2| \, dx} \\
&& \le 2\|n(|I|)\|_{L^2}(\|\tilde{u}(|I|)\|_{L^2}+\|w(|I|)\|_{L^2})\|w(|I|)\|_{L^{\infty}} \\
&& \le c\|n(|I|)\|_{L^2}(\|\tilde{u}(|I|)\|_{L^2}+\|w(|I|)\|_{L^2})\|w(|I|)\|_{L^2}^{\frac{1}{2}} \|w_x(|I|)\|_{L^2}^{\frac{1}{2}} \\
&& \le cN^{1-s}(1+N^{-\frac{3}{2}+\frac{3}{4}s+})N^{\frac{1}{2}(-\frac{3}{2}+\frac{3}{4}s)+}N^{\frac{1}{2}(\frac{3}{2}-2s)+} \\
&& \le cN^{1-s}N^{-\frac{5}{8}s+}
\end{eqnarray*}
It is easy to see that the decisive bound is the one for the first term. Thus the increment of the energy is bounded by
$$ c_3 N^{1-s}N^{\frac{3}{2}-2s+} $$
where $ c_3 = c_3(\bar{c},M) $.\\
Thus the condition which ensures uniform control of the energy of $(\tilde{u},\tilde{n}),(\tilde{\tilde{u}},\tilde{\tilde{n}}),\dots $ is the following:
\begin{equation}
\label{104a}
c_3TN^{4(1-s)+}N^{1-s}N^{\frac{3}{2}-2s+} < \bar{c} N^{2(1-s)} 
\end{equation}
where $ c_3 = c_3(2\bar{c},2M) $ (recall that by (\ref{50a}) initially the energy is bounded by $ \bar{c}N^{2(1-s)} $ ). \\
This is fulfilled for $N$ sufficiently large provided
$$ 4(1-s)+(1-s)+\frac{3}{2}-2s < 2(1-s) \quad \Longleftrightarrow \quad s > \frac{9}{10} $$
So, here is the point where the decisive bound on $s$ appears. The uniform control of the energy implies by (\ref{49a}),(\ref{49b}) uniform control of the $L^2$-norm of $ (\tilde{u}_x,\tilde{n},A^{-1/2}\tilde{n}_t),(\tilde{\tilde{u}}_x,\tilde{\tilde{n}},A^{-1/2}\tilde{\tilde{n}}_t),\dots $
\begin{theorem}
Let $1>s> 9/10$. The Zakharov system (\ref{0.1}),(\ref{0.1a}),(\ref{0.2}) with data $(u_0,n_0,n_1) \in H^s({\bf R})\times L^2({\bf R}) \times \dot{H}^{-1}({\bf R})$ is globally well-posed. More precisely for any $T>0$ there exists a unique solution
\begin{equation}
\label{5.5}
(u,n,n_t) \in X^{s,\frac{1}{2}+\epsilon_1}[0,T] \times X^{0,\frac{1}{2}+\epsilon_2}[0,T] \times \dot{X}^{-1,\frac{1}{2}+\epsilon_2}[0,T]
\end{equation}
for $\epsilon_1,\epsilon_2>0$ small enough. This solution satisfies
\begin{equation}
\label{5.6}
(u,n,n_t) \in C^0([0,T],H^s({\bf R})\times L^2({\bf R}) \times \dot{H}^{-1}({\bf R}))
\end{equation}
\end{theorem}
{\bf Proof:} On any of the intervals $I$ of the preceding considerations we have by (\ref{76}),(\ref{3.16}),(\ref{3.12}) (+ interpolation) $ \tilde{u} \in X^{s,\frac{1}{2}+\epsilon_1}(I),$ $\tilde{n}_{\pm} \in X^{0,\frac{1}{2}+\epsilon_2}(I)$ and by (\ref{4.10}),(\ref{4.12}) $ v \in X^{s,\frac{1}{2}+\epsilon_1}(I),$ $m_{\pm} \in X^{0,\frac{1}{2}+\epsilon_2}(I)$. This gives (\ref{5.5}) by (\ref{2.3b}),(\ref{2.3c}). Uniqueness in this class was proven in \cite{GTV} already. (\ref{5.6}) follows immediately from (\ref{5.5}).

It is not difficult to give bounds on the growth of the solutions now. It is elementary to show that the most restrictive bound on $N$ comes from condition (\ref{104a}) in the whole range $9/10 < s<1$, namely
\begin{equation}
\label{5.6a}
N > cT^{\frac{1}{5s-\frac{9}{2}}+}
\end{equation}
According to the construction of the solution above we have the following structure
$$ u(t)=\tilde{u}(t)+e^{it\partial_x^2}u_{02} + w(t)= e^{it\partial_x^2}u_0 + r(t) $$
where
$$ r(t) = \tilde{u}(t)-e^{it\partial_x^2}u_{01} + w(t) $$
on $I$ first, but then also on $[0,T]$.\\
We have shown that
\begin{equation}
\|r(t)\|_{H^{1,2}({\bf R})} \le cN^{1-s}
\label{5.7}
\end{equation}
(remark that $\|e^{it\partial_x^2}u_{01}\|_{H^{1,2}({\bf R})} \le cN^{1-s}$). Choosing $N$ according to (\ref{5.6a}) gives the following bound for $0\le t \le T$:
$$ \|r(t)\|_{H^{1,2}({\bf R})} \le cT^{\frac{1-s}{5s-\frac{9}{2}}+} $$
Similarly
$$ \|n_{\pm}(t)\|_{L^2({\bf R})} \le cN^{1-s} \le cT^{\frac{1-s}{5s-\frac{9}{2}}+}  $$
Thus we have shown
\begin{theorem}
The solution of the preceding theorem fulfills for $t \ge 0$:
$$ u(t)= e^{it\partial_x^2}u_0 + r(t) $$
with
$$ \|r(t)\|_{H^{1,2}({\bf R})} \le 
c \left(1+|t|^{\frac{1-s}{5s-\frac{9}{2}}+}\right) $$
and
$$ \|n(t)\|_{L^2({\bf R})}+ \|n_t(t)\|_{\dot{H}^{-1,2}({\bf R})} \le 
c \left(1+|t|^{\frac{1-s}{5s-\frac{9}{2}}+}\right)$$
\end{theorem}

\end{document}